\documentclass[12pt]{article}
\usepackage{latexsym,amssymb,amsmath}
\usepackage{amsthm, amstext}
\usepackage{array, amsfonts, mathrsfs}
\usepackage{graphicx}

\newcommand{\range}[2]{{#1:#2}}
\newcommand{\md}{\mathrm{d}}

\newtheorem{convention}{Convention}

\numberwithin{equation}{section}

\begin{document}
\date{}

\author{M.I.Belishev\thanks {St.Petersburg Department of
                 the Steklov Mathematical Institute, St.Petersburg State University, Russia;
                 belishev@pdmi.ras.ru. Supported by the grants
                 RFBR 14-01-00535À and SPbGU 6.38.670.2013.}, I.B.Ivanov\thanks{St.Petersburg State University,
                 Institute of Physics, St.Petersburg Nuclear Physics Institute, Theoretical Physics Division, Russia;
                 contact@ivisoft.org}, I.V.Kubyshkin\thanks{St.Petersburg State University, Institute of
                 Physics, Russia; kub@geo.phys.spbu.ru}, V.S.Semenov\thanks{St.Petersburg State University, Institute of
                 Physics, Russia; sem@geo.phys.spbu.ru. Supported by grant of the RF President for support of leading scientific schools 2836.2014.5.}}

\title{Numerical testing in determination of sound speed from a part of boundary by the
BC-method}
\maketitle

 \begin{abstract}
We present the results of numerical testing on determination of
the sound speed $c$ in the acoustic equation $u_{tt}-c^2\Delta
u=0$ by the {\it boundary control method}. The inverse data is a
response operator (a hyperbolic Dirichlet-to-Neumann map) given on
controls, which are supported on {\it a part} of the boundary. The
speed is determined in the subdomain covered by acoustic rays,
which are emanated from the points of this part orthogonally to
the boundary. The determination is {\it time-optimal}: the longer
is the observation time, the larger is the subdomain, in which $c$
is recovered. The numerical results are preceded with brief
exposition of the relevant variant of the BC-method.
 \end{abstract}



\section{Introduction}\label{sec:intro}

\subsection{About the method}\label{sec: About the method}
The {\it boundary control method} (BCM) is an approach to inverse
problems based on their relations with control and system theory
\cite{BIP97,BIP07,Isbook}. It is a rigorously justified
mathematical method of synthetic character: Riemannian geometry,
asymptotic methods in PDE, functional analysis and operator theory
are in the use. Beginning on its foundation in 1986 \cite{BDAN87},
there was a question whether such a purely theoretical method is
available for numerical implementation. The first affirmative
results were obtained by V.B.Filippov in two-dimensional problem
of the density $\rho=c^{-2}$ determination via the spectral
inverse data \cite{BRF}. Later on, an algorithm based on the
spectral variant of the BCM was elaborated and tested by
S.A.Ivanov and V.Yu.Gotlib in \cite{BGICOCV,BIP07}.

A dynamical variant of the BCM deals with time-domain inverse data
that is a response operator (hyperbolic Dirichlet-to-Neumann map).
It provides {\it time-optimal} reconstruction: the longer is the
observation time, the bigger is the subdomain, in which the
parameters are recovered. It is the feature, which makes this
variant most relevant for possible applications to acoustics and
geophysics. The corresponding algorithm was elaborated and tested
by V.Yu.Gotlib in \cite{BGotJIIPP}. It recovers the density in a
near-boundary layer from the data given on the {\it whole
boundary}.

Time-optimal determination of density via the spectral and
time-domain inverse data given on a {\it part of boundary} is
proposed in \cite{BIP97}. The procedure uses singular harmonic
functions; its spectral variant was realized numerically (see
\cite{BIP97}, section 7.7). In \cite{BHow02} and \cite{CUBO}, its
dynamical variant was modified to make it more prospective for
applications in geophysics, the modification being based on
geometrical optics.

In beginning of 2000's, L.Pestov proposed a version of the BCM,
which determines some intrinsic bilinear forms containing
parameters under reconstruction via the inverse data and, then,
recovers the parameters from the forms. This version is not
time-optimal but, on expense of big enough observation time,
provides more stable numerical algorithms. The results of the
collaboration, which develops this approach in the I.Kant Baltic
Federal University (Kaliningrad, Russia), are presented in \cite{P-12,P-13,PKBolg}.

Recently, L.Oksanen applied the BCM for numerical reconstruction
of the obstacle \cite{Oks}.

There also exists a time optimal and {\it data optimal} approach
by V.Romanov \cite{Rom} but it is not implemented and tested yet.
Another (not optimal) direct reconstruction methods, which are
numerically (and experimentally) tested, see in \cite{BeiKl, KSh,
KSh_1}.

\subsection{Inverse problem}\label{sec:Inverse problem}
The goal of our work is to elaborate the BC-algorithm for
time-optimal determination of the sound speed via the time-domain
inverse data given at a part of boundary, and test it in numerical
experiment.
\smallskip

$\bullet$\,\, Let $\Omega \subset {\mathbb R}^n$ be a (possibly,
unbounded) domain with the boundary $\Gamma$. We deal with a
dynamical system
 \begin{align}
\label{acoustic1} &u_{tt} - c^2\Delta u = 0 && {\text {in}}\,\,\Omega\times(0,T)\\
\label{acoustic2} &u|_{t=0} = u_t|_{t=0} = 0 && {\text {in}}\,\,{\overline \Omega} \\
\label{acoustic3} &u = f && {\text {on}}\,\,\Gamma\times[0,T]\,,
 \end{align}
where $c=c(x)$ is a smooth enough positive function ({\it speed of
sound}), $f$ is a {\it boundary control}, $u=u^f(x,t)$ is a
solution ({\it wave}). With the system one associates a {\it
response operator}
 \begin{equation}\label{R^T}
R^T:\, f \mapsto u^f_\nu\big|_{\Gamma \times [0,T]}\,,
 \end{equation}
where $(...)_\nu$ is a derivative with respect to the outward
normal $\nu$ on $\Gamma$. In a general form, the inverse problem
is to answer the question: To what extent does the response
operator determine the sound speed into the domain? Also, the
determination procedures are of principal interest.

System (\ref{acoustic1})--(\ref{acoustic3}) is hyperbolic and, as
such, obeys the {\it finiteness of the domains of influence}
(FDI). It describes the waves propagating with finite speed $c$,
and the relevant setup of the inverse problem must take this
property into account. Such a setup is given below, after
geometric preliminaries.
\smallskip

$\bullet$\,\,The sound speed induces a {\it travel time metric}
$d\tau^2=c^{-2}|dx|^2$ (shortly, $c$-metric) and the corresponding
distance $\tau(x,y)$ in $\Omega$. For a subset $A \subset
\overline \Omega$, by
 \begin{equation*}
\Omega^\xi_A:=\{x\in \overline \Omega\,|\,\,\tau(x,A)<\xi\}
 \end{equation*}
we denote its $c$-metric neighborhood of radius $\xi$.
\smallskip

By $r^\xi_\gamma$ we denote a geodesic in $c$-metric ({\it ray}),
which is emanated from $\gamma \in \Gamma$ into $\Omega$ in
direction $-\nu$, and is of the $c$-length $\xi$. Let $\sigma
\subset \Gamma$ be a part of the boundary. A set
 \begin{equation*}
B^T_\sigma:=\bigcup_{\gamma \in \sigma}r^T_\gamma\subset \overline
{\Omega^T_\sigma}
 \end{equation*}
is called a ray {\it tube}. On Fig 1a,b, the neighborhood
$\Omega^T_\sigma$ and tube $B^T_\sigma$ are contoured by the
closed lines $\{1,2,3,4,5,6,1\}$ and $\{5,6,2,3,5\}$ respectively
($B^T_\sigma$ is shaded).

\begin{figure}[ht]
\centering
\includegraphics[width=5in,height=2in]{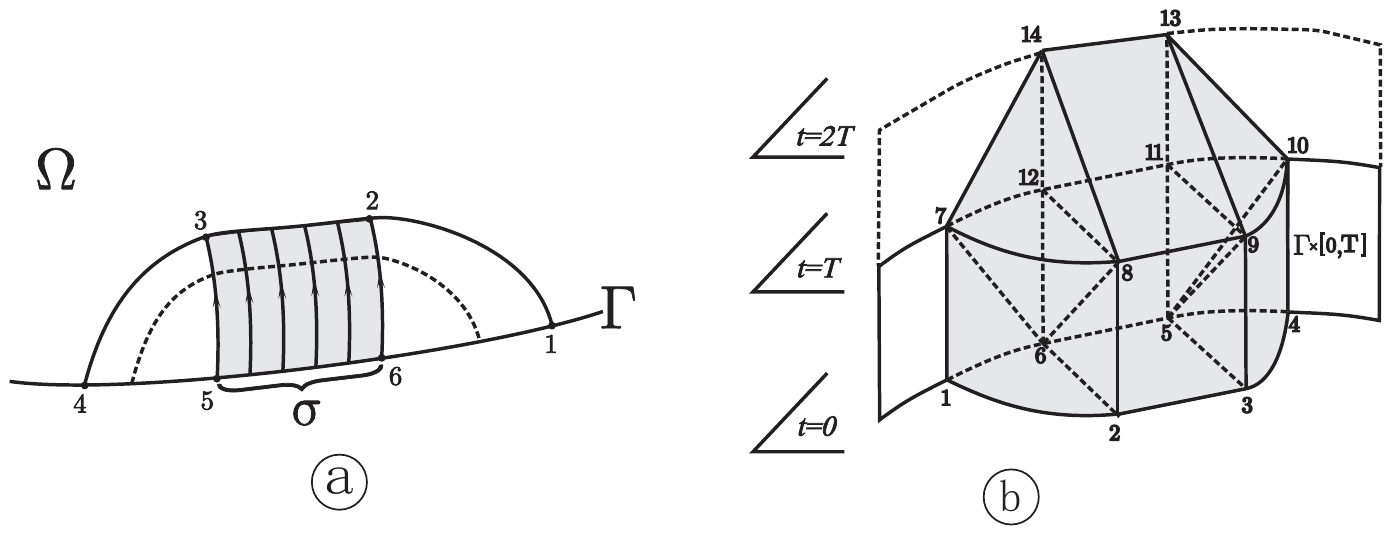}
\caption{Tube $B^T$ and domain $D^{2T}_\sigma$} \label{fig:Tube}
\end{figure}

If $T$ is small enough then the ray field is regular in the tube.
Let $T_\sigma$ be the infimum of $T$'s, for which such a
regularity does occur.
 \begin{convention}\label{convention sigma disk}
In what follows, unless otherwise specified, we assume that
$\sigma$ is diffeomorphic to a disk $\{p \in {\mathbb
R}^{n-1}\,|\,\,|p|\leq 1\}$ and $T<T_\sigma$. Such a case is said
to be regular.
 \end{convention}

The part $\sigma$ determines the space-time domains
 \begin{align*}
& D^{2T}_\sigma:=\{(x,t)|\,x \in \Omega^T_\sigma,\,
0<\tau(x,\sigma)<2T-t\} \quad {\rm and}\\
& E^T_\sigma:=\{(x,t)|\,x \in \Omega^T_\sigma,\, 0\le
t<\tau(x,\sigma)\},
 \end{align*}
and the space-time surfaces
 \begin{align*}
\Theta^{T}_\sigma:=\overline{D^{2T}_\sigma}\cap\{\Gamma\times[0,T]\}\,,
\quad
\Theta^{2T}_\sigma:=\overline{D^{2T}_\sigma}\cap\{\Gamma\times[0,2T]\}.
 \end{align*}
All of them are mapped by the projection $(x,t)\mapsto x$ to
$\overline{\Omega^T_\sigma}$. Domain $D^{2T}_\sigma$ is shown on
Fig 1.b (shaded). Domain $E^{T}_\sigma \subset D^{2T}_\sigma$ lies
under the surface $\{(x,t)\,|\,\,t=\tau(x,\sigma)\}$, which
consists of three parts countered by the closed lines
$\{6,7,8,6\}$, $\{5,6,8,9,5\}$, and $\{5,9,10,5\}$. The surfaces
$\Theta^{2T}_\sigma$ and $\Theta^{T}_\sigma$ are countered by the
lines $\{1,6,5,4,10,13,14,7,1\}$ and $\{1,6,5,4,10,11,12,7,1\}$
respectively.
\smallskip

If $c<c_*=\rm const$ holds in $\overline \Omega$ then for the sets
 \begin{equation}\label{sigma^T}
\sigma^\xi:=\{\gamma \in \Gamma\,|\,\,\tau(\gamma,\sigma)\leq
\xi\} \quad {\rm and}\quad\sigma^\xi_*:=\{\gamma \in
\Gamma\,|\,\,{\rm dist}_{{\mathbb R}^n}(\gamma,\sigma)\leq
c_*\xi\}
 \end{equation}
one has $\sigma^\xi \subset \sigma^\xi_*$, and the relations
 \begin{equation}\label{embeddings sigma^xi in sigma^T}
\Theta^{T}_\sigma\,=\,\sigma^T\times[0,T]\,\subset\, \sigma^T_*
\times [0,T]
 \end{equation}
are valid.
\smallskip

$\bullet$\,\,Assign a control $f$ to a class ${\cal
F}^{2T}_\sigma$ if ${\rm supp\,}f \subset \sigma \times [0,2T]$,
i.e., it acts from $\sigma$ during the time interval $0\leq t\leq
2T$. Owing to the FDI, an extension of system
(\ref{acoustic1})--(\ref{acoustic3}) of the form
\begin{align}
\label{acoustic1*} &u_{tt} - c^2\Delta u = 0 && {\text {in}}\,\,D^{2T}_\sigma\\
\label{acoustic2*} &u = 0 && {\text {in}}\,\,E^{T}_\sigma\\
\label{acoustic3*} &u = f \in {\cal F}^{2T}_\sigma
 \end{align}
turns out to be a well-posed problem, its solution $u^f$ being
determined by the values of the speed $c$ in the subdomain
$\Omega^T_\sigma$ (does not depend on
$c|_{\Omega\setminus\Omega^T_\sigma}$). The same is valid for the
response operator
 \begin{equation}\label{R^2Tsigma}
R^{2T}_\sigma: \,f \mapsto u^f_\nu\big|_{\Theta^{2T}_\sigma}
 \end{equation}
associated with this problem: it is also determined by
$c|_{\Omega^T_\sigma}$.

By the latter, the relevant setup of the inverse problem is: {\it
for a fixed $T>0$, given the operator $R^{2T}_\sigma$ determine
the speed $c$ in $\Omega^T_\sigma$.}
\smallskip

The use of the doubled time $2T$ is quite natural by kinematic
reasons. The subdomain $\Omega^T_\sigma$ is prospected with waves
initiated at $\sigma$. To search the whole $\Omega^T_\sigma$, the
waves have to fill it (that takes $T$ time units) and return back
to the boundary (for the same time $T$) to be detected by the
external observer, which implements measurements at $\Gamma$.
 \begin{convention}\label{convention extend by zero}
The operator $R^{2T}_\sigma$ is introduced so that, for the times
$0\leq t \leq T$ the images $R^{2T}_\sigma f$ are defined on the
set $\Theta^T_\sigma$ only. For convenience of further
formulations, we put $R^{2T}_\sigma f \big|_{0\leq t \leq T}$ to
be extended from $\Theta^T_\sigma$ to $\Gamma \times [0,T]$ by
zero.
 \end{convention}

\subsection{Results and comments}\label{Results and comments}
$\bullet$\,\,Let $\sigma \subset \Gamma$ and $T>0$ be given. Our
{\it a priori assumptions} are that $T<T_\sigma$ (i.e., we deal
with the regular case) and the sound speed upper bound $c_*$ is
known. Under these assumptions, we propose a procedure, which
recovers the speed $c$ in the tube $B^T_\sigma$ via the operator
$R^{2T}_\sigma$. Then, we demonstrate the results of numerical
testing of the algorithm based on this procedure.
\smallskip

$\bullet$\,\, In fact, the procedure utilizes not the complete
operator $R^{2T}_\sigma$ but some information, which it
determines. Namely, as will be seen, to recover
$c\big|_{\Omega^T_\sigma}$, it suffices for the external observer
to possess the following options:
 \begin{enumerate}
 \item for any $f,g \in {\cal F}^{2T}_\sigma$ obeying the oddness
condition
 $$
f(\cdot,T)=-f(\cdot,2T-t),\quad g(\cdot,T)=-g(\cdot,2T-t)\,,\qquad
0\leq t \leq 2T\,,
 $$
one can compute the integral
 \begin{align}\label{form phi}
I^T_\sigma[f,g]:=\int_{\sigma\times [0,2T]}u^{Jf}_\nu(\gamma,t)
g(\gamma,t)\,d\Gamma dt=(R^{2T}_\sigma Jf,g)_{{\cal
F}^{2T}_\sigma}\,,
 \end{align}
where $J: {\cal F}^{2T}_\sigma \to {\cal F}^{2T}_\sigma$ is an
integration: $(Jf)(\cdot,t):=\int_0^t f(\cdot, s)\,ds$.

 \item for any odd $f \in {\cal F}^{2T}_\sigma$, one can detect
$u^f_\nu\big|_{\sigma^T_*\times [0,T]}=R^{2T}_\sigma
f\big|_{\sigma^T_*\times [0,T]}$, i.e., implement the measurements
on $\sigma^T_*$ (but not on the whole $\Gamma$!) during the time
in\-ter\-val $[0,T]$ (but not $[0,2T]$!)
 \end{enumerate}
\smallskip

$\bullet$\,\,In principle, the proposed procedure is identical to
the versions \cite{BHow02} and \cite{CUBO}. Therefore, its
exposition is short: we omit some proofs and derivations,
referring the reader to the mentioned papers for detail. In the
mean time, here we deal with more refined (rigorously
time-optimal) data that is the operator $R^{2T}_\sigma$, in
contrast to \cite{BHow02} and \cite{CUBO}, where the operator
$R^{2T}$ corresponding to system
(\ref{acoustic1})--(\ref{acoustic3}) with the final time $t=2T$,
is used as the inverse data.
\smallskip

$\bullet$\,\,One of the features and advantages of the BCM is that
it reduces nonlinear inverse problems to linear ones. In
particular, the main fragment of the algorithm, which recovers
$c$, is the solving a big-size {\it linear} algebraic system. The
matrix of the system is of the form
$\{I^T_\sigma[f_i,f_j]\}_{i,j=1}^N$ for a rich enough set of
controls $f_i$. As a consequence of the strong ill-posedness of
the above stated inverse problem, this system also turns out to be
ill posed but the linearity enables one to apply standard
regularization devices.

\section{Geometry}\label{sec:geom}

\subsection{$c$-metric}\label{sec:c-metric}
$\bullet$\,\, Let $\Omega \subset {\mathbb R}^n$ ($n\geq 2$) be a
domain with the $C^2$-smooth boundary $\Gamma$.  A {\it sound
speed} is a function $c \in C^2(\overline \Omega)$ provided $c>0$.
If $\Omega$ is unbounded, we assume $c \leq c_*={\rm const}$.
\smallskip

The sound speed determines a {\it $c$-metric} in $\Omega$ with the
length element $d\tau^2=c^{-2}dl^2$ and the distance
 \begin{align*}
\tau(x,y) := \text{inf} \int_x^y\,\frac{dl}{c}\,,
 \end{align*}
where $dl$ is the Euclidean length element, and the infimum is
taken over smooth curves, which lie in $\Omega$ and connect $x$
with $y$. In dynamics, the value $\tau(x,y)$ is a travel time
needed for a wave initiated at $x$ to reach $y$.
\smallskip

$\bullet$\,\,\, Let $\sigma\subset\Gamma$; a function
 \begin{align*}
\tau_\sigma(x) := \inf_{y\in\sigma}\,\tau(y,x)\,, \qquad x \in
\Omega
 \end{align*}
is called an {\it eikonal}. Its value is the travel time from
$\sigma$ to $x$. A set
 \begin{align*}
\Omega^\xi_\sigma:=\left\{x\in\Omega \mid \tau_\sigma(x) <
\xi\right\} \qquad (\xi>0)
 \end{align*}
is a $c$-metric neighborhood of $\sigma$ of radius $\xi$. In
dynamics, the waves initiated on $\sigma$ at the moment $t=0$,
fill up the subdomain $\Omega^\xi_\sigma$ at $t=\xi$. The filled
domains are bounded by the eikonal level sets
\begin{align*}
\Gamma^\xi_\sigma:=\left\{x\in\Omega \mid \tau_\sigma(x) =
\xi\right\}
\end{align*}
(the surfaces $c$-equidistant to $\sigma$: see the dotted line on
Fig 1a), which play the role of the forward fronts of waves
propagating from $\sigma$ into $\Omega$.

\subsection{Ray coordinates}\label{Ray coordinates}
$\bullet$\,\,\,Fix a point $\gamma\in\sigma$. Let
$x(\gamma,\xi)\in \Omega$ be the endpoint of the $c$-metric
geodesic ({\it ray}) $r^\xi_\gamma$ starting from $\gamma$
orthogonally to $\Gamma$ and parametrized by its $c$-length $\xi$.
Also, put $x(\gamma,0)\equiv \gamma$.

For $T>0$, the rays starting from $\sigma$, cover a {\it tube}
 \begin{align*}
B^T_\sigma\,=\,\bigcup_{\gamma\in\sigma}r^T_\gamma=\bigcup_{(\gamma,
\xi)\in\sigma\times[0,T]}
x(\gamma,\xi)\,\subset\,\overline{\Omega^T_\sigma}\,.
 \end{align*}
In the regular case, $B^T_\sigma$ is diffeomorphic to the set
$$\Sigma^T_\sigma\,:=\,\sigma \times [0,T]$$  via the map
$\Sigma^T_\sigma \ni(\gamma, \xi) \mapsto x(\gamma, \xi) \in
B^T_\sigma$ (on Fig 1.b, $\Sigma^T_\sigma$ is countered by
$\{6,5,11,12,5\}$ ). This enables one to regard a pair
$(\gamma,\xi)$ as the {\it ray coordinates} of the point
$x(\gamma,\xi)\in B^T_\sigma$.
\smallskip

$\bullet$\,\,Let $\pi^i$ be the Cartesian coordinate functions:
$\pi^i(x):=x^i$ for $x=\{x^i\}_{i=1}^n \in {\mathbb R}^n$. The map
 \begin{equation}\label{ray to Cart}
(\gamma, \xi)\mapsto \{\pi^i\left(x(\gamma,
\xi)\right)\}_{i=1}^n\,, \qquad (\gamma,\xi)\in \Sigma^T_\sigma
 \end{equation}
realizes the passage from the ray coordinates to Cartesian ones.

Fix a $\gamma \in \sigma$. The equality
 \begin{equation}\label{c=sqrt dot x^2}
c\left(x(\gamma,\xi)\right)\,=\,\left\lbrace\sum \limits_{i=1}^n
\left[\frac{d}{d\xi}\,\pi^i\left(x(\gamma,
\xi)\right)\right]^2\right \rbrace^\frac{1}{2}, \qquad
0<\xi<T
 \end{equation}
represents $c$ on the ray $r^T_\gamma$. Varying $\gamma$, we get
the sound speed representation in the whole tube
$\Sigma^T_\sigma$.

\subsection{Images}\label{def images}
$\bullet$\,\, Fix a point $\gamma\in\sigma$, choose a small
$\varepsilon>0$, and define the surfaces
 $$
\sigma_\varepsilon(\gamma, \xi):=\{x(\gamma',\xi) \in
B^T_\sigma\,|\,\,\tau(\gamma',\gamma)< \varepsilon\}\,, \qquad
0\leq \xi< T\,.
 $$
A function
 \begin{align*}
J(\gamma,\xi) := \lim_{\varepsilon\to
0}\,\frac{|\sigma_\varepsilon(\gamma,
\xi)|}{|\sigma_\varepsilon(\gamma, 0)|}\,, \qquad (\gamma, \xi)\in
\Sigma^T_\sigma,
 \end{align*}
where $|...|$ is a surface area in ${\mathbb R}^n$, is said to be
a {\it ray spreading} at the point $x(\gamma,\xi)$.

In the regular case, the coefficients
 $$
\varkappa(\gamma,\xi):=\frac{J(\gamma,\xi)}{c(x(\gamma,\xi))}\quad
\text{and} \quad
\beta(\gamma,\xi):=\left[\varkappa(\gamma,0)\varkappa(\gamma,\xi)\right]^{\frac{1}{2}}\,,
 $$
which enter in the well-known geometrical optics relations (see,
e.g., \cite{BB,Ik}), are the smooth functions on
$\Sigma^T_\sigma$.
\smallskip

$\bullet$\,\, Let $y$ be a function on $B^T_\sigma$; a function
$\tilde y$ of the form
 \begin{align*}
\tilde y(\gamma,\xi):=\beta(\gamma,\xi)\,y(x(\gamma,\xi)),\qquad
(\gamma, \xi)\in \Sigma^T_\sigma
 \end{align*}
is called an {\it image} of $y$. For the function $\pi^0(x)\equiv
1$, one has $\tilde \pi^0 = \beta$.

In terms of images, relations (\ref{ray to Cart}) and (\ref{c=sqrt
dot x^2}) take the form of the representations
 \begin{align}
\notag &(\gamma, \xi)\mapsto \left\{\frac{\tilde\pi^i(\gamma,
\xi)}{\tilde\pi^0(\gamma, \xi)}\right\}_{i=1}^n=x(\gamma, \xi),
\quad c\left(x(\gamma,\xi)\right)\,=\,\left\lbrace\sum
\limits_{i=1}^n
\left[\frac{d}{d\xi}\left(\frac{\tilde\pi^i(\gamma,
\xi)}{\tilde\pi^0(\gamma, \xi)}\right)\right]^2\right
\rbrace^\frac{1}{2},\\
\label{BASIC} & (\gamma,\xi)\in \Sigma^T_\sigma\,,
 \end{align}
which will be used for determination of $c$ in the inverse
problem.

\section{Dynamics}\label{sec:bc}
In section \ref{sec:bc}, the regularity condition $T<T_\sigma$ is
cancelled, and $T>0$ is arbitrary. However, for the sake of
simplicity, we keep $\sigma$ to be diffeomorphic to a disk. All
the functions, spaces, operators, etc are {\it real}. We denote
$\Sigma^s_\sigma:=\sigma \times [0,s]$.

\subsection{Spaces and operators}\label{Spaces and operators}
Denote the dynamical system associated with problem
(\ref{acoustic1})--(\ref{acoustic3}) by $\alpha^T$. In what
follows, we deal with its subsystem corresponding to controls
acting from $\sigma$. We consider it as a separate system, denote
by $\alpha^T_\sigma$, and endow with standard control theory
attributes: spaces and operators. All of them are determined by
$c\big|_{\Omega^T_\sigma}$.
\smallskip

$\bullet$\,\,The space of boundary controls $\mathcal{F}^T_\sigma
:= L_2(\Sigma^T_\sigma)$ with the inner product
 \begin{align*}
\left(f,g\right)_{\mathcal{F}^T_\sigma} :=
\int_{\Sigma^T_\sigma}f(\gamma,t)\,g(\gamma,t)\,d\Gamma dt
 \end{align*}
($d\Gamma$ is the Euclidean surface element on the boundary) is
called an {\it outer space} of system $\alpha^T_\sigma$. It
contains an increasing family of subspaces
 \begin{align*}
\mathcal{F}^{T,\xi}_\sigma := \left\{f\in \mathcal{F}^T_\sigma
\mid {\rm supp\,}f \subset \sigma\times[T-\xi,T]\right\},\qquad
0\le \xi\le T
 \end{align*}
(${\cal F}^{T,0}_\sigma=\{0\},\,\,{\cal F}^{T,T}_\sigma={\cal
F}^{T}_\sigma$) formed by the delayed controls acting from
$\sigma$. Here, $T-\xi$ is the value of delay, $\xi$ is an action
time.
\smallskip

$\bullet$\,\,The space $\mathcal{H}^T_\sigma :=
L_2(\Omega^T_\sigma; c^{-2}dx)$ with the inner product
 \begin{align*}
\left(y,w\right)_{\mathcal{H}^T_\sigma} := \int_{\Omega^T_\sigma}
y(x)\,w(x)\,\frac{dx}{c^2(x)}
 \end{align*}
is said to be an {\it inner} space of the system. It contains a
family of subspaces
 \begin{align*}
\mathcal{H}^{\xi}_\sigma \,:=\, \left\{y\in\mathcal{H}^T_\sigma
\mid \text{supp}\,y\subseteq
\overline{\Omega^\xi_\sigma}\right\},\quad 0\le \xi\le T
 \end{align*}
(${\cal H}^{0}_\sigma:=\{0\}$), which increase as $\sigma$ extends
and/or $\xi$ grows.
\smallskip

$\bullet$\,\, In the system $\alpha^T_\sigma$, an
`input$\rightarrow$state' correspondence is described by a {\it
con\-t\-rol operator} $W^T:\mathcal{F}^T_\sigma\rightarrow
\mathcal{H}^T_\sigma$,
$$W^T f\,:=\, u^f(\cdot,T)\,,$$
where $u^f$ is a solution to (\ref{acoustic1})--(\ref{acoustic3}).
Operator $W^T$ is bounded \cite{BIP97}.

Since the waves governed by the equation (\ref{acoustic1})
propagate with the finite speed $c$, for controls acting from
$\sigma$ one has
 \begin{align}\label{supp u^f}
{\rm supp\,}u^f(\cdot,\xi) \subset \overline{\Omega^\xi_\sigma}\,,
\qquad 0\leq \xi \leq T \,.
 \end{align}
As is easy to recognize, (\ref{supp u^f}) is equivalent to the
embedding
 \begin{align}\label{Embedding}
W^T\mathcal{F}^{T,\xi}_\sigma\subset \mathcal{H}^\xi_\sigma,
\qquad 0\leq \xi \leq T\,.
  \end{align}

$\bullet$\,\,Recall that $\nu$ is the outward normal to $\Gamma$,
and the sets $\sigma^\xi$ are defined in (\ref{sigma^T}).  Denote
$\Sigma^T:=\Gamma \times [0,T]$

An `input$\to$output' correspondence is realized by the {\it
response operator} $R^T:\mathcal{F}^T_\sigma\to L_2(\Sigma^T;
d\Gamma dt)$,
 \begin{equation*}
R^{T} f \,:=\, u^f_\nu\big|_{\Sigma^T}
 \end{equation*}
defined on the set ${\rm Dom\,}R^T=\{f \in
H^1(\Sigma^T_\sigma)\,\big|\,\,f\big|_{\partial \sigma \times
[0,T]}=0,\,\,f\big|_{t=0} = 0\}$, where $H^1(...)$ is the Sobolev
class and $\partial \sigma$ is the boundary of $\sigma$ in
$\Gamma$). Relation (\ref{supp u^f}) implies
 \begin{align*}
{\rm supp\,}u^f_\nu\, \subset \,\{(\gamma,\xi)\,|\,\,\gamma \in
\sigma^\xi,\,\,\,0\leq \xi \leq T \}\subset \sigma^T\times
[0,T]\,.
 \end{align*}
By the latter, for controls $f \in {\cal F}^T_\sigma$, one has
\begin{align}\label{supp R^Tf*}
{\rm supp\,}R^Tf\,\subset\,\sigma^T\times
[0,T]\overset{(\ref{embeddings sigma^xi in sigma^T})}\subset
\,\sigma^T_*\times [0,T]\,.
 \end{align}

One more (extended) response operator $R^{2T}_\sigma: {\cal
F}^{2T}_\sigma \to L_2(\Theta^{2T}_\sigma; d\Gamma dt)$ is
\begin{equation*}
R^{2T}_\sigma f\,:=\,u^f_\nu\big|_{\Theta^{2T}_\sigma}\,,
 \end{equation*}
where $u^f$ is a solution to extended problem
(\ref{acoustic1*})--(\ref{acoustic3*}). It is defined on ${\rm
Dom\,}R^{2T}_\sigma=\{f \in
H^1(\Sigma^{2T}_\sigma)\,\big|\,\,f\big|_{\partial \sigma \times
[0,2T]}=0,\,\,f\big|_{t=0} = 0\}$. As was noted in
\ref{sec:Inverse problem}, $R^{2T}_\sigma$ is determined by the
values of the sound speed $c$ in the subdomain $\Omega^T_\sigma$.
Therefore, it is reasonable to regard it as an intrinsic object of
system $\alpha^T_\sigma$ (but not $\alpha^{2T}_\sigma$!).
\smallskip

Let the controls $f\in {\rm Dom\,}R^T$ in (\ref{acoustic3}) and
$\check f\in{\rm Dom\,}R^{2T}_\sigma$ in (\ref{acoustic3*}) be
such that $f=\check f\big|_{0\leq t \leq T}$. Then, the solutions
to problems (\ref{acoustic1})--(\ref{acoustic3}) and
(\ref{acoustic1*})--(\ref{acoustic3*}) also coincide for the same
times:
 \begin{align*}
u^f=u^{\check f} \qquad {\rm in}\,\,\Omega^T_\sigma \times
[0,T]\,.
 \end{align*}
As a consequence, passing to the normal derivatives on $\Gamma$,
one gets
 \begin{equation}\label{R^T and R^2Tsigma}
R^T f\,=\,R^{2T}_\sigma \check f \quad {\rm
on}\,\Sigma^T\cap\Theta^{2T}_\sigma\,.
 \end{equation}

$\bullet$\,\, A {\it connecting operator} of the system is
$C^T:\mathcal{F}^T_\sigma\rightarrow \mathcal{F}^T_\sigma$,
 \begin{equation*}
C^T\,:=\,(W^T)^\ast W^T\,.
 \end{equation*}
The definition implies
\begin{align}\label{Cop}
\left(u^f(\cdot,T), u^g(\cdot,T)\right)_{\mathcal{H}^T_\sigma} =
\left(W^T f, W^Tg\right)_{\mathcal{H}^T_\sigma} = \left(C^T f,
g\right)_{\mathcal{F}^T_\sigma}\,,
\end{align}
i.e., $C^T$ connects the Hilbert metrics of the outer and inner
spaces.

A significant fact is that the connecting operator is determined
by the response operator in a simple explicit way. Namely, the
representation
\begin{equation}\label{C T via R 2T}
 C^T\,=\,2^{-1}(S^T)^* R^{2T}_\sigma J S^T
\end{equation}
is valid, where the map $S^T: {\cal F}^T_\sigma \to {\cal
F}^{2T}_\sigma$ extends the controls from $\sigma \times [0,T]$ to
$\sigma \times [0,2T]$ by oddness with respect to $t=T$:
 \begin{align}
\left(S^Tf\right)(\cdot,t)\,:=\,
\begin{cases}
\,\,\,\,f(\cdot,t), &0\le t<T\\-f(\cdot,2T-t), &T\le t \le 2T
\end{cases}
 \end{align}
and $J: {\cal F}^{2T}_\sigma\to {\cal F}^{2T}_\sigma$ is an
integration:\,\,$(Jf)(\cdot,t)=\int_0^t f(\cdot,s)\,ds $ (see
\cite{BIP97}--\cite{CUBO}). Note that $f\in{\rm Dom\,}R^T$ implies
$S^Tf\in{\rm Dom\,}R^{2T}_\sigma$.

As a consequence, we get
 \begin{align}
\notag & \left(C^T f,g\right)_{{\cal
F}^{T}_\sigma}\overset{(\ref{Cop})}=\left(u^f(\cdot,T),
u^g(\cdot,T)\right)_{{\cal
 H}^T_\sigma}\overset{(\ref{C T via R 2T})}=
2^{-1}\left(R^{2T}_\sigma JS^T f,S^Tg\right)_{{\cal
 F}^{2T}_\sigma}=\\
\label{!!!} &  \overset{(\ref{form
phi})}=\,2^{-1}I^T_\sigma[S^Tf,S^Tg]
 \end{align}
for arbitrary $f \in {\rm Dom\,}R^{T}$ and $g \in {\cal
F}^{T}_\sigma$.

\subsection{Wave bases}\label{sec:wavbas}
For the BCM, the fact of crucial character is that the embedding
(\ref{Embedding}) is dense: the equality
 \begin{align}\label{LBC}
{\rm
clos\,}W^T\mathcal{F}^{T,\xi}_\sigma\,=\,\mathcal{H}^\xi_\sigma,
\qquad 0\le \xi\le T
 \end{align}
(the closure in $\mathcal{H}^T_\sigma$) is valid and interpreted
as a {\it local boundary controllability} of system
(\ref{acoustic1})--(\ref{acoustic3}). It shows that the waves
constitute rich enough sets in the subdomains which they fill up.
In particular, by this property, any square-summable function
supported in $\Omega^T_\sigma$ can be approximated (with any
precision) by a wave $u^f(\cdot, T)$ owing to proper choice of the
control $f$ acting from $\sigma$
\cite{BHow02,BIP07,BBlag99,BGotJIIPP}.
\smallskip

$\bullet$\,\,An important consequence of controllability is
existence of wave bases.

Fix a $\xi \in (0,T]$. Let a linearly independent system of
controls $\{f^\xi_k\}_{k=1}^\infty$ be {\it complete} in the
subspace ${\cal F}^{T,\xi}_\sigma$, i.e. the relation $\vee
\{f^\xi_k\}_{k=1}^\infty= {\cal F}^{T,\xi}_\sigma$ holds, where
$\vee$ is a closure of the linear span (in the relevant norm). By
(\ref{LBC}), the system of waves
 $$
u^\xi_k:=u^{f^\xi_k}(\cdot,T)=W^Tf^\xi_k
 $$
turns out to be complete in ${\cal H}^\xi_\sigma$, i.e., one has
  $
\vee \{u^\xi_k\}_{k=1}^\infty= {\cal H}^\xi_\sigma\,.
  $

If $T$ is such that $\Omega \setminus \Omega^T_\sigma
\not=\emptyset$, i.e., the waves moving from $\sigma$ do not cover
the whole $\Omega$, then the control operator is injective
\cite{ABI} (in particular, this holds for $T<T_\sigma$). In this
case, $W^T$ preserves the linear independence, and
$\{u^\xi_k\}_{k=1}^\infty$ turns out to be a linearly independent
complete system in ${\cal H}^\xi_\sigma$.
 \begin{convention}\label{convention wave basis}
By this, we deal with this case and say $\{u^\xi_k\}_{k=1}^\infty$
to be a {\it wave basis} in the subspace ${\cal H}^\xi_\sigma$.
Also, everywhere, system $\{f^\xi_k\}_{k=1}^\infty$ producing the
wave basis, is chosen so that all $f^\xi_k \in {\rm Dom\,}R^T$.
 \end{convention}
\noindent As a consequence, the Gramm marices
  $${\cal G}^\xi_N\,:=\,\{(u^\xi_i, u^\xi_j)_{{\cal H}^T_\sigma}\}_{i,j=1}^N, \qquad N=1,2,\dots$$
are nonsingular and invertible, whereas their entries can be
represented via the controls:
 \begin{equation}\label{Gramm G}
({\cal G}^\xi_N)_{ij} =(C^Tf^\xi_i, f^\xi_j)_{{\cal
F}^T_\sigma}\overset{(\ref{!!!})}= 2^{-1}(R^{2T}_\sigma JS^T
f^\xi_i,S^Tf^\xi_j)_{{\cal
 F}^{2T}_\sigma}\,.
 \end{equation}
\smallskip

$\bullet$\,\,In the BCM, wave bases are used for finding the
projections of functions on the domains filled with waves.

Fix a positive $\xi \leq T$. Let $P^\xi_\sigma$ be the
(orthogonal) projector in ${\cal H}^T_\sigma$ onto ${\cal
H}^\xi_\sigma$. Such a projector cuts off functions:
 \begin{equation*}
P^\xi_\sigma y\,=\,\begin{cases}
 y & {\rm in}\,\,\Omega^\xi_\sigma\\
 0 & {\rm in}\,\,\Omega^T_\sigma \setminus\Omega^\xi_\sigma
                   \end{cases}\,.
 \end{equation*}
As an element of the subspace ${\cal H}^\xi_\sigma$, this
projection can be represented via the wave basis:
 \begin{equation}\label{Proj on WB}
P^\xi_\sigma y\,=\,\lim \limits_{N\to \infty}P^\xi_{\sigma,N}=\lim
\limits_{N\to \infty}\sum \limits_{k=1}^Nc^\xi_{k,N}u^\xi_k\,,
 \end{equation}
where $P^\xi_{\sigma,N}$ projects in ${\cal H}^T_\sigma$ onto the
span $\vee\{u^\xi_k\}_{k=1}^N$, and the column of coefficients
$\{c^\xi_{k,N}\}_{k=1}^N=:C^\xi_N$ is determined via the column
$\{(y, u^\xi_k)_{{\cal H}^T_\sigma}\}_{k=1}^N=:B^\xi_N$ through
the Gramm matrix by
 $$
C^\xi_N\,=\,\left[{\cal G}^\xi_N\right]^{-1}B^\xi_N\,.
 $$
The limit is understood in the sense of the norm convergence in
${\cal H}^T_\sigma$.

\subsection{Dual system}
$\bullet$\,\,Denote $ K^{T}_\sigma:=\{(x,t)\,|\,\,0<
\tau(x,\sigma)<t<T\}$.

\begin{figure}[ht]
\centering
\includegraphics[width=3.3in,height=2.2in]{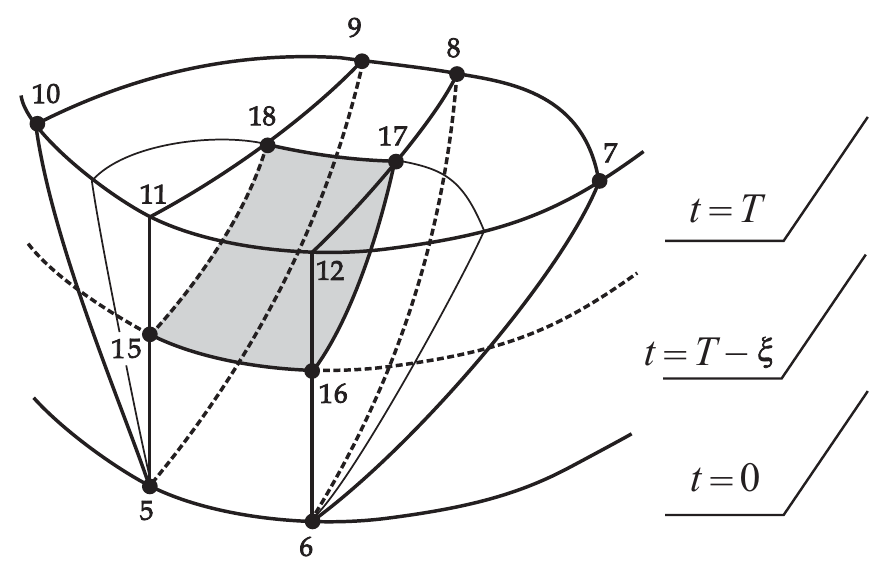}
\caption{Domain $K^{T}_\sigma$} \label{fig:Dual}
\end{figure}

A dynamical system associated with the problem
 \begin{align}
\label{dual1} & v_{tt} - c^2\Delta v = 0 && {\text {in}}\,\,\,K^{T}_\sigma\\
\label{dual2} & v|_{t=T} = 0,\,\,\,v_t|_{t=T} = y && {\text {in}}\,\,\,\overline {\Omega^T_\sigma} \\
\label{dual3} & v = 0 && {\text {on}}\,\,\,\sigma^T\!\times\![0,T]
 \end{align}
is called {\it dual} to system $\alpha^T_\sigma$; by $v=v^y(x,t)$
we denote its solution. Owing to the FDI, such a problem turns out
to be well possed for any $y \in {\cal H}^T_\sigma$. Its
peculiarity is that the Cauchy data are assigned to the {\it
final} moment $t=T$, so that the problem is solved in reversed
time.

The solutions to the original and dual problems obey the {\it
duality relation}: for any $f \in {\cal F}^T_\sigma$ and $y \in
{\cal H}^T_\sigma$, the equality
\begin{align}\label{duality}
\left(u^f(\cdot,T), y\right)_{{\cal H}^T_\sigma}\, = \,\left(f,
v^y_\nu\right)_{\mathcal{F}^T_\sigma}
\end{align}
is valid \cite{BIP97,BHow02,BIP07,BBlag99}.
\smallskip

$\bullet$\,\,With the dual system one associates an {\it
observation operator} $O^T: {\cal H}^T_\sigma \to {\cal
F}^T_\sigma$,
 $$
O^Ty\,:=\,v^y_\nu \big|_{\Sigma^T_\sigma}\,.
 $$
Writing (\ref{duality}) in the form $\left(W^Tf,y\right)_{{\cal
H}^T_\sigma}\, = \,\left(f, O^T y\right)_{\mathcal{F}^T_\sigma}$,
we get an  operator equality $O^T\,=\,(W^T)^*$. Hence, the
definition of $C^T$ implies
\begin{equation}\label{C^T=OW}
C^T\,=\,O^T W^T\,.
 \end{equation}

\subsection{Projections of harmonic functions}\label{sec:Projections of harmonic functions}
$\bullet$\,\,Assume that $y$ in (\ref{dual2}) is harmonic: $y=a
\in {\cal H}^T_\sigma$ obeys $\Delta a=0$ in $\Omega^T_\sigma$ and
is continuously differentiable up to $\sigma^T \subset \partial
\Omega^T_\sigma$. A simple integration by parts in (\ref{duality})
leads to
 \begin{align}
\notag & \left(a,u^f(\cdot,T)\right)_{\mathcal{H}^T_\sigma}
=\left(O^T a, f\right)_{\mathcal{F}^T_\sigma}=\\
\label{(a,u^f)} & =\int_{\sigma^T
\times[0,T]}(T-t)\left[a(\gamma)(R^Tf)(\gamma,t)-a_\nu(\gamma)f(\gamma,t)\right]\,d\Gamma
dt
 \end{align}
(see \cite{BIP97,BHow02,BBlag99,BGotJIIPP}).

Assume that $\{f^\xi_k\}_{k=1}^\infty\subset{\cal
F}^{T,\xi}_\sigma$ is chosen in accordance with Convention
\ref{convention wave basis} and produces the wave basis
$\{u^\xi_k\}_{k=1}^\infty\subset{\cal H}^{\xi}_\sigma$. Then,
representation (\ref{Proj on WB}) takes the form
 \begin{equation}\label{P^xia=sum}
P^\xi_\sigma a\,=\,\lim \limits_{N\to \infty}\sum
\limits_{k=1}^Nc^\xi_{k,N} u^\xi_k\,,
 \end{equation}
where $C^\xi_N=\{c^\xi_{k,N}\}_{k=1}^N$ satisfies the linear
system
 \begin{equation}\label{Eqn GC=B}
{\cal G}^\xi_N C^\xi_N\,=\,B^\xi_N
 \end{equation}
with the Green matrix
  \begin{align}
\notag &{\cal G}^\xi_N\,\overset{(\ref{Gramm
G})}=\,\,\left\{2^{-1}(R^{2T}_\sigma JS^T
f^\xi_i,S^Tf^\xi_j)_{{\cal
 F}^{2T}_\sigma}\right\}_{i,j=1}^N\,=\\
\label{Gramm G1} & =\left\{2^{-1}\int_{\sigma\times
[0,2T]}(R^{2T}_\sigma JS^T
f^\xi_i)(\gamma,t)\,(S^Tf^\xi_j)(\gamma,t)\,d\Gamma
dt\right\}_{i,j=1}^N
 \end{align}
and the right-hand side
 \begin{align*}
\notag &  B^\xi_N=\{(a, u^\xi_k)_{\cal H}\}_{k=1}^N\,, \quad (a,
u^\xi_k)_{\cal H}\overset{(\ref{(a,u^f)})}=\\
 &
=\,\int_{\sigma^T
\times[0,T]}(T-t)\left[a(\gamma)(R^Tf^\xi_k)(\gamma,t)-
a_\nu(\gamma)f^\xi_k(\gamma,t)\right]\,d\Gamma dt\,.
 \end{align*}
With regard to (\ref{supp R^Tf*}),(\ref{R^T and R^2Tsigma}), and
Convention \ref{convention extend by zero}, the latter can be
written in the form
 \begin{align}\label{B}
B^\xi_N=\left\{\int_{\sigma^T_*\times
[0,T]}(T-t)\left[a(\gamma)(R^{2T}_\sigma S^T f^\xi_k)(\gamma,t)-
a_\nu(\gamma)f^\xi_k(\gamma,t)\right]\,d\Gamma dt\right\}_{k=1}^N
 \end{align}
determined by $R^{2T}_\sigma$ and, thus, relevant for the further
use.
\smallskip

$\bullet$\,\, Fix a positive $\xi<T$. The operator
 $$
P^\xi_{\sigma\bot}:=P^T_\sigma - P^\xi_\sigma
 $$
is the projector in $\cal H$ onto the subspace ${\cal
H}^T_\sigma\ominus{\cal H}^\xi_\sigma$; it cuts off functions on
the subdomain $\Omega^T_\sigma \setminus \Omega^\xi_\sigma$.

Choose systems $\{f^T_k\}_{k=1}$ and $\{f^\xi_k\}_{k=1}$, which
are linearly independent and complete in ${\cal F}^T_\sigma$ and
${\cal F}^{T,\xi}_\sigma$ respectively. Applying the (bounded)
observation operator to (\ref{P^xia=sum}), with regard to
$O^Tu^f=O^T W^T f\overset{(\ref{C^T=OW})}=C^Tf$, we obtain
  \begin{equation*}
O^TP^T_\sigma a\,=\,\lim \limits_{N\to \infty}\sum
\limits_{k=1}^Nc^T_{k,N}C^T f^T_k, \quad O^TP^\xi_\sigma
a\,=\,\lim \limits_{N\to \infty}\sum \limits_{k=1}^Nc^\xi_{k,N}C^T
f^\xi_k\,.
 \end{equation*}
Subtracting,  we arrive at the representation
 \begin{equation}\label{O(P^T-P^xi)a=sum}
O^TP^\xi_{\sigma \bot} a\,=\,\lim \limits_{N\to \infty}\sum
\limits_{k=1}^N \left[c^T_{k,N}C^T f^T_k - c^\xi_{k,N}C^T
f^\xi_k\right]\,.
 \end{equation}
For the future application to the inverse problem, a crucial fact
is that its right-hand side {\it is determined by the response
operator}. Indeed, if $R^{2T}_\sigma$ is given, one can
 \begin{enumerate}
\item choose the complete linearly independent systems
$\{f^T_k\}^\infty_{k=1} \subset {\cal F}^T_\sigma$ and
$\{f^\xi_k\}^\infty_{k=1}\subset {\cal F}^{T, \xi}_\sigma$; then,
com\-pose the Gramm matrices ${\cal G}^T_N$, ${\cal G}^\xi_N$ by
(\ref{Gramm G1}) and columns $B^T_N$, $B^\xi_N$ by (\ref{B})

\item solving system (\ref{Eqn GC=B}) with respect to $C^T_N$,
$C^\xi_N$, find the coefficients $c^T_{k,N},\,c^\xi_{k,N}$

\item determine $C^T$ by (\ref{C T via R 2T}), compose the sum in
(\ref{O(P^T-P^xi)a=sum}) and, extending $N$, pass to the limit.
 \end{enumerate}

\subsection{Amplitude formula}\label{AF}
In what follows, we deal with the regular case $T<T_\sigma$.
\smallskip

$\bullet$\,\, Fix a positive $\xi<T$; let $y$ be a smooth function
in $\Omega$. Return to the dual system
(\ref{dual1})--(\ref{dual3}) and put
$$v_t \big|_{t=T}\,=\,P^\xi_{\sigma \bot} y\,=:\,y^\xi_\bot\,$$ in Cauchy data (\ref{dual2}).
Such a $y^\xi_\bot$ is of two specific features:
\smallskip

\noindent($i$)\,\,\,it vanishes in $\Omega^\xi_\sigma$, so that
${\rm supp}\,y^\xi_\bot$ is separated from $\sigma$ by the
$c$-distance $\xi$. Therefore, by the finiteness of the wave
propagation speed, $v^{y^\xi_\bot}$ vanishes in the space-time
domain $\{(x,t)\in K^T_\sigma\,|\,\,t>(T-\xi)+\tau(x, \sigma)\}$
and, in particular, one has
\begin{equation}\label{separated from sigma}
v^{y^\xi_\bot}_\nu(\cdot, t)\big|_\sigma \,=\,0 \qquad {\rm
for}\,\,T-\xi<t\leq T\,.
 \end{equation}

\noindent($ii$)\,\,$y^\xi_\bot$ is discontinuous: generically, it
has jumps at the equidistant surfaces $\Gamma^T_\sigma$ and
$\Gamma^\xi_\sigma$. In particular, at the points
$x(\gamma,\xi)\in B^T_\sigma \cap \Gamma^\xi_\sigma$, the value
({\it amplitude}) of the jump is
 \begin{equation}\label{jump1}
y^\xi_\bot(x(\gamma,\xi+0))\,=\,y(x(\gamma,\xi))\,.
 \end{equation}

$\bullet$\,\, In hyperbolic equations theory, the well-known fact
is that discontinuous Cauchy data initiate discontinuous
solutions, the discontinuities propagating along characteristics.
In our case of the wave equation (\ref{dual1}), the jumps of
$v_t\big|_{t=T}=y^\xi_\bot$ induce the jumps of $v^{y^\xi_\bot}_t$
in $K^T_\sigma$. In particular,  there is a jump on the
characteristic surface
 $\{(x,t)\in \overline{K^T_\sigma}\,|\,\,t=T-\xi+\tau(x,
\sigma)\}$ including its smooth part
 $
S^{T,\xi}_\sigma:=\{(x,t)\in \overline{K^T_\sigma}\,|\,\,x \in
B^T_\sigma\}
 $
(on Fig 2, contoured by $\{15,16,17,18,15\}$). The jumps of
$v^{y^\xi_\bot}_t$ on $S^{T,\xi}_\sigma$ and of
$v^{y^\xi_\bot}_\nu$ on the cross-section $S^{T,\xi}_\sigma
\cap\Sigma^T_\sigma=\{(\gamma, T-\xi)\,|\,\,\gamma \in \sigma\}$
(the line $\{15,16\}$) can be found by standard geometrical optics
devices. For the latter jump, a simple analysis provides
 \begin{equation*}
v^{y^\xi_\bot}_\nu(\gamma,t)\big|^{t=T-\xi+0}_{t=T-\xi-0}\,=\,-
\beta(\gamma,\xi)y(x(\gamma,\xi))\,, \qquad \gamma \in \sigma
 \end{equation*}
(see, e.g., \cite{BB},\cite{Ik},\cite{BBlag99}). By
(\ref{separated from sigma}), we have
$v^{y^\xi_\bot}_\nu(\gamma,t)\big|^{t=T-\xi+0}=0$ that leads to
 \begin{equation}\label{jump2}
v^{y^\xi_\bot}_\nu(\gamma,T-\xi-0)\,=\,\beta(\gamma,\xi)y(x(\gamma,\xi))\,,
\qquad \gamma \in \sigma\,.
 \end{equation}
Comparing (\ref{jump1}) with (\ref{jump2}), one can recall the
well-known physical principle: jumps propagate along rays (here, a
ray is $r^\xi_\gamma=\{x(\gamma,s)\,|\,\, 0 \leq s \leq \xi\}$)
with the speed $c$, the ratio of the jump amplitudes at the input
and output of the ray (here, at $x(\gamma, \xi)$ and
$x(\gamma,0)$) depending on the ray spreading.
\smallskip

$\bullet$\,\, Recalling the definitions of images and observation
operator, one can write (\ref{jump2}) in the form
 \begin{equation}\label{AmplForm}
(O^T P^\xi_{\sigma \bot} y)(\gamma,\xi)\,=\,\tilde y(\gamma,
\xi)\,, \qquad (\gamma, \xi) \in \Sigma^T_\sigma\,.
 \end{equation}
It is the so-called {\it amplitude formula} (AF), which plays a
central role in solving inverse problems by the BCM \cite{BIP97,
BIP07, BBlag99}. It represents the image of function in the form
of collection of jumps, which pass through the medium, absorb
information on the medium structure, and are detected by the
external observer at the boundary.
\smallskip

$\bullet$\,\, Now, let $y=a$ be a harmonic function. Combining
(\ref{O(P^T-P^xi)a=sum}) with (\ref{AmplForm}), we arrive at the
key relation
 \begin{align}
\notag &\tilde a(\gamma,\xi)\,=\,\lim \limits_{t \to
T-\xi-0}\left\lbrace\lim \limits_{N\to \infty}\left(\sum
\limits_{k=1}^N \left[c^T_{k,N}\,C^T f^T_k - c^\xi_{k,N}\,C^T
f^\xi_k\right]\right)(\gamma,t)\right\rbrace\,,\\
\label{key rel} & (\gamma, \xi) \in \Sigma^T_\sigma\,.
 \end{align}
As was noted at the end of section \ref{sec:Projections of
harmonic functions}, to find its right-hand side, it suffices to
know the response operator. In particular, since the coordinate
functions are {\it harmonic}, applying (\ref{key rel}) to
$a=\pi^i, \,\,i=0, \dots, n$ one can recover their images $\tilde
\pi^i$ via $R^{2T}_\sigma$.

\section{Determination of speed}\label{Determination of speed}
\subsection{Procedure}\label{sec:procedure}
To solve the inverse problem, we just summarize our considerations
in the form of the following procedure. Recall that the role of
the procedure input data is played by operator $R^{2T}_\sigma$.
\smallskip

\noindent{\bf Step 1.}\,\, Fix a $\xi<T$. Applying the procedure
$1.-3.$ described at the end of section \ref{sec:Projections of
harmonic functions}, find the right-hand side of (\ref{key rel})
for $a=\pi^0, \pi^1, \dots, \pi^n$ and, thus, get the images
$\tilde \pi^i(\gamma, \xi)$ for $\gamma \in \sigma$.
\smallskip

\noindent{\bf Step 2.}\,\,Varying $\xi$, find $\tilde \pi^i$ on
$\Sigma^T_\sigma$. Then, recover the map $\Sigma^T_\sigma \ni
(\gamma,\xi) \mapsto x(\gamma,\xi)\in {\mathbb R}^n$ by the first
representation in (\ref{BASIC}). The image of the map is
$B^T_\sigma$, so that the ray tube is recovered in $\Omega$.
\smallskip

\noindent{\bf Step 3.}\,\, Differentiating with respect to $\xi$,
find $c$ by the second representation in (\ref{BASIC}). The pairs
$\{x(\gamma,\xi),c(x(\gamma,\xi))\,|\,\,(\gamma,\xi) \in
\Sigma^T_\sigma\}$ constitute the graph of $c$ in $B^T_\sigma$.

Thus, the sound speed in  the tube is determined. The following is
some comments and remarks.
\smallskip

$\bullet$\,\,For applications in geophysics, by the obvious
reasons, it is desirable to minimize the part of the boundary, on
which the external observer has to implement measurements. As is
seen from (\ref{B}), our procedure requires observations not only
on $\sigma$ but on $\Gamma \backslash \sigma$, whereas the
knowledge of the bound $c_*$ just enables the observer to restrict
measurements on $\sigma^T_*$. In principle, one can avoid the
observations on $\Gamma \backslash \sigma$ by the use of the
artificial coordinates instead of the Cartesian $\pi^i$. Namely,
one can choose the harmonic functions $a^1, \dots, a^n$ obeying
$a^i\big|_{\Gamma \backslash \sigma}=0$, which separate points of
the tube $B^T_\sigma$ at least locally. By this choice, in
(\ref{B}) one gets
$\int_{\sigma^T_*\!\times\![0,T]}=\int_{\sigma\!\times\![0,T]}$.
Therefore, possessing the values of $R^{2T}_\sigma f$ on
${\sigma\!\times\![0,2T]}$ (but not on the whole
$\Theta^{2T}_\sigma$!), one can recover the ima\-ges $\tilde a^i$
via the amplitude formula and use them for identifying the points
of $B^T_\sigma$ in $\Omega$.  Thereafter, one recovers
$c\big|_{B^T_\sigma}$. However, it is not clear, whether this plan
can provide workable numerical algorithms.
\smallskip

$\bullet$\,\, As was mentioned in \ref{sec: About the method}, the
procedure \cite{BIP97} (sections 7.6, 7.7) enables one to
determine $c\big|_{\Omega^T_\sigma}$ from observations on
${\sigma\!\times\![0,2T]}$ {\it only}, and, thus, provides the
strongest uniqueness result. However, its numerical implementation
in the case of the time-domain inverse data seems to be rather
problematic.

\subsection{Numerical testing} \label{sec:num res}
\subsubsection*{Preparation of tests}
$\bullet$\,\, We take
 \begin{align*}
& \Omega:=\{(x^1,x^2)\in {\mathbb R}^2\,|\,\,x^2\le 0\}\,,\quad
\Gamma:=\{(x^1,0)\in {\mathbb R}^2\,|\,\,-\infty < x^1 <
\infty\}\,,\\
& \sigma:=\{x\in\Gamma\,|\,\,-L\le x^1\le L\}\,,
\end{align*}
and consider a few concrete examples of the {\it density}
$\rho=c^{-2}$ in $\Omega$.

\noindent$\bullet$\,\,We choose an appropriate finite system of
controls $f_k$ supported on $\sigma\times [0,T]$, the system being
{\it the same for all examples} except Test $1(b)$ where results
with another spatial basis are presented for comparison. It is
well adapted to constructing the systems of delayed controls: for
intermediate $\xi=\xi_l$, the shifts $f_k^{\xi_l}(\cdot,
t)=f_k(\cdot, t-(T-\xi_l))$ are in use. This enables one to reduce
considerably the computational resources.

\noindent$\bullet$\,\, At each of the examples, we solve
numerically the forward problems
(\ref{acoustic1})--(\ref{acoustic3}) with the final moments $t=T$
and $t=2T$ for the controls $f^{\xi_l}_i$ and ${JS^Tf^{\xi_l}_i}$
respectively. These problems are solved by the use of a
semi-discrete central-upwind third order accurate numerical scheme
with WENO reconstruction suggested in \cite{kurganov2001}. As a
result, we get the functions $u^{f^{\xi_l}_i}_\nu=R^Tf^{\xi_l}_i=
(R^{2T}_\sigma S^Tf^{\xi_l}_i)\big|_{0 \le t \le T}$ and
$u^{JS^Tf^{\xi_l}_i}_\nu=R^{2T}_\sigma JS^Tf^{\xi_l}_i$ entering
in (\ref{B}) and (\ref{Gramm G1}).

\subsubsection*{Controls}\label{sec:basis}
The BCM uses a system of boundary controls $f_1,\,f_2,\,\dots$,
which belong to the Sobolev class:
\begin{equation}
\nonumber \left\{f_k\in H^1(\Gamma\times [0,T]) \mid
f_k(\gamma,t)|_{t=0} = 0\right\}
\end{equation}
and constitute a basis in $L_2(\Gamma\times [0,T])$. We construct
such a system from the products of elements of spatial and
temporal bases, $f_k(\gamma,t) = \phi_l(\gamma)\,\psi_m(t)$, $k=l
+ m N_\gamma$, where $l=\range{0}{N_\gamma-1}$,
$m=\range{0}{N_t-1}$, and the basis dimension is $N = N_\gamma
N_t$.

\begin{figure}[ht]
\includegraphics[width=0.478\textwidth]{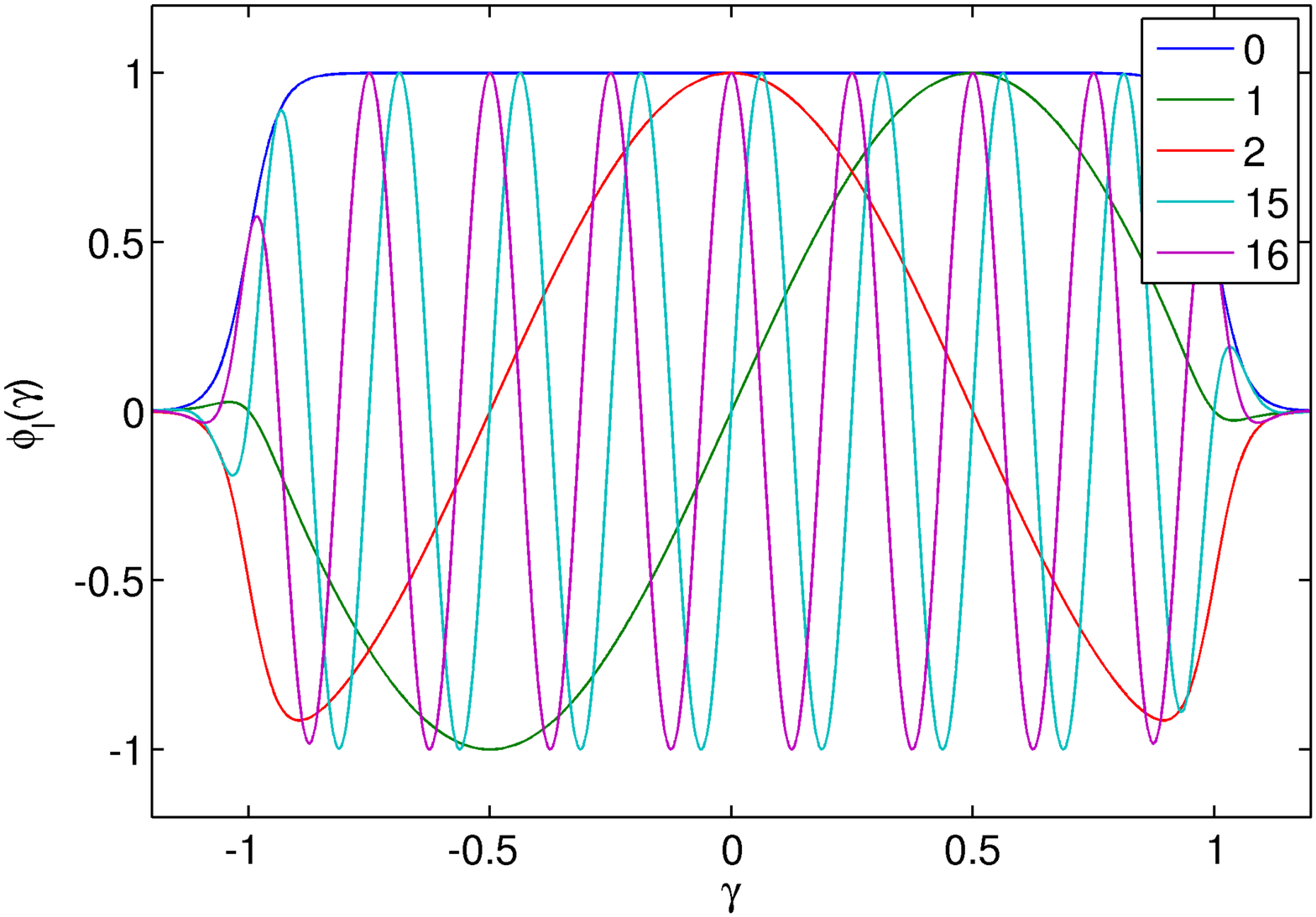}\hfill
\includegraphics[width=0.473\textwidth]{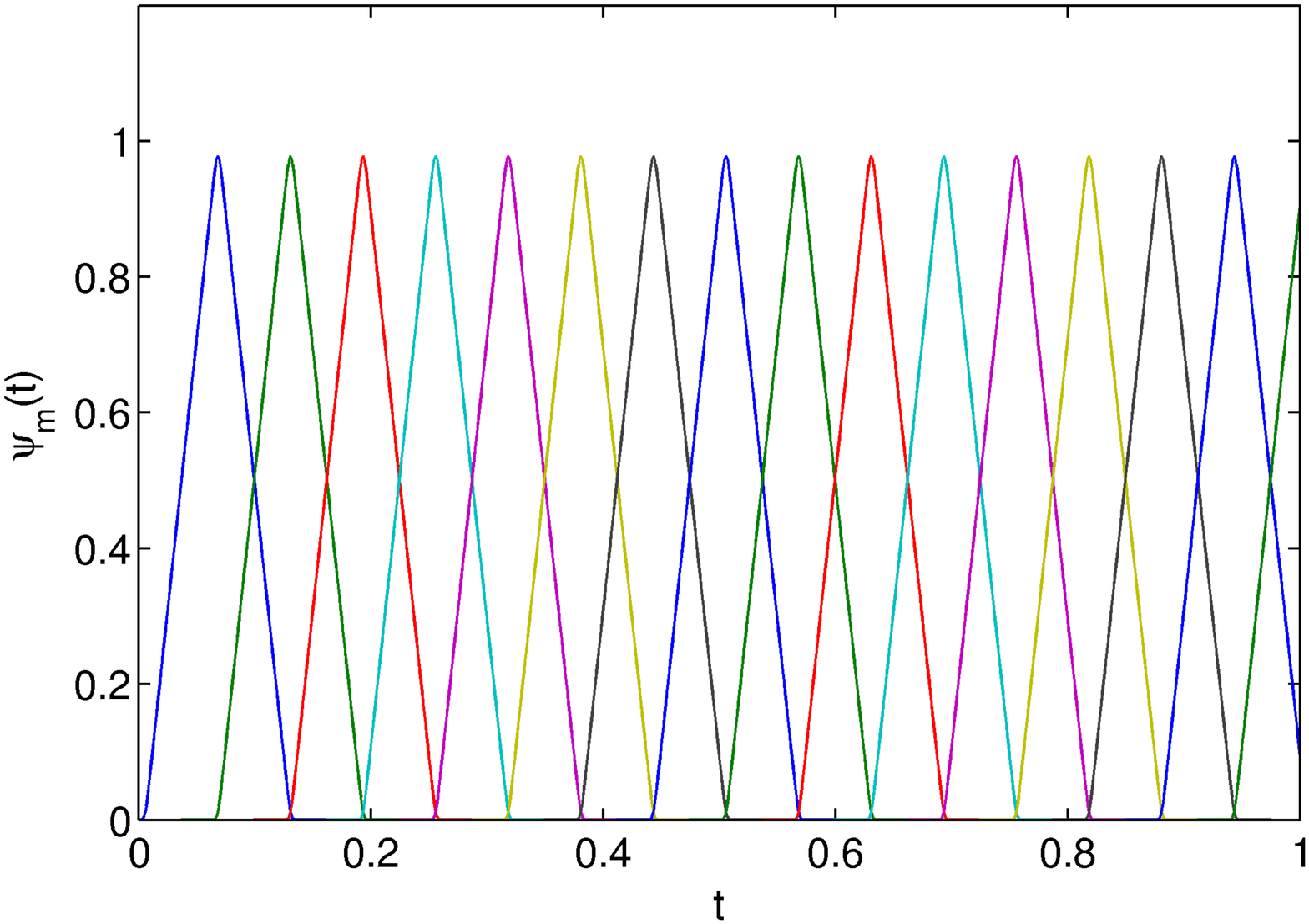}
\caption{Basis of boundary controls: spatial functions $\phi_l(\gamma)$ with $l=\range{0}{16}$ and $s=1/32$ (left)
and temporal functions $\psi_m(t)$ with $m=\range{0}{15}$ and $T = 1$, $N_t = 16$, $d=\Delta/64$ (right).}
\label{fig:basis}
\end{figure}

In the case of the half-plane, we can keep under control only a
part of the boundary and thus have to use {\it localized} basis
functions. The simplest and good choice is a conventional
trigonometric basis reduced to the interval $[-1, 1]$ by an
exponential cutoff multiplier $\eta(\gamma) = 1 / (1 +
\exp\left(\gamma/s\right))$ with a cutoff scale $s$, so that
\begin{equation}
\label{spatial} \phi_l(\gamma) = \eta(\gamma-1)\,
\eta(-\gamma-1)\,\cos\left[\pi\left(\frac{l}{2} + \lfloor\frac{l +
1}{2}\rfloor (\gamma - 1)\right)\right]\,,
\end{equation}
where $\lfloor \cdot \rfloor$ is the integer part. The spatial
basis functions are shown in the left panel of
Fig.~\ref{fig:basis}.

The temporal basis is constructed from the shifts of a {\it
tent-like} function,
\begin{equation}
\label{temporal}
\theta(t) = \frac{d}{\Delta}\,
\left(1 - \exp\left[-\frac{\Delta}{d}\right]\right)^{-1}
\ln\left[
\frac
{
\cosh\left[\frac{2\Delta-t}{2 d}\right]
\cosh\left[\frac{t}{2 d}\right]
}
{
\cosh^2\left[\frac{\Delta-t}{2 d}\right]
}
\right]
,
\end{equation}
so that $\psi_m(t) = \theta(t - m\Delta - \delta)$, where $\Delta
= T / N_t$, $d$ is a smoothing parameter (when $d\rightarrow 0$
the function $\theta(t)$ gets a triangular shape), and $\delta$ is
an offset to ensure a negligible value of $\theta(0)$. Such a
shift-invariant basis (shown in the right panel of
Fig.~\ref{fig:basis}) considerably reduces computational resources
needed for the BCM-reconstruction.

\subsubsection*{Regularization}\label{sec:regularization}
$\bullet$\,\, In the course of determination of $c$ by the
procedure {\bf Step 1-3}, we use the above-prepared data for
computing the entries of $B^{\xi_l}_N$ and ${\cal G}^{\xi_l}_N$ in
(\ref{B}) and (\ref{Gramm G1}). Then the system (\ref{Eqn GC=B})
is solved for $a=\pi^0, \pi^1, \pi^2$, and the solutions $C^{\xi_l}_N$
are calculated by standard LAPACK routines.

Solving system (\ref{Eqn GC=B}), we have to apply a regularization
procedure since the condition number of the Green matrix ${\cal
G}^{\xi_l}_N$ rapidly grows as its size $N$ increases, see Figure
\ref{fig:case2:cond}. Because of unavoidable errors in matrix
elements and right hand sides, the expansion coefficients
$C^{\xi_l}_N$ also contain errors amplified by ill-conditioned
matrix. We use Tikhonov's regularization to reduce fake
oscillations caused by errors in expansion coefficients. The value
of regularization parameter is selected to satisfy a desired
tolerance for residual of the linear system.

\smallskip

\noindent$\bullet$\,\, One more operation, which produces
unavoidable errors, is computation of the double limit in
(\ref{key rel}). The origin of the errors is the following.

In (\ref{P^xia=sum}), the projection $P^\xi_\sigma a$ is a
piece-wise smooth function in $\Omega$, which has a jump at the
surface $\Gamma^\xi_\sigma$. Therefore, the convergence of the
sums in the right hand side not uniform near $\Gamma^\xi_\sigma$,
and the Gibbs oscillations do occur in the summation process.
These oscillations are transferred to the amplitude formula
(\ref{key rel}) and considerably complicate the determination of
jump at $t=T-\xi$, whereas this determination is a crucial point
of the algorithm.

To damp this negative effect we apply the following procedure. The
basis functions have finite resolution of the order of
spatial-temporal scales of the highest harmonic. All scales below
the minimum ones are unreachable, therefore we {\it average} the
result of expansions (\ref{key rel}) over that minimum scales by
convolution with some kernel $K(\gamma, \xi)$,
\begin{align}
\label{convol} \langle g\rangle(\gamma,t) =
\int\limits_{-\infty}^{+\infty}\md
t'\int\limits_{-\infty}^{+\infty}\md\gamma' K(\gamma - \gamma', t
- t')\,g(\gamma', t').
\end{align}
In our implementation the kernel $K(\gamma - \gamma', t - t')$ is
a product of conventional {\it Gaussian kernels} both for spatial
and temporal variables. Such a procedure efficiently removes the
Gibbs oscillations and, in fact, accelerates convergence of the
expansions. The values of standard deviations in the Gaussian
kernels should match the minimum spatial and temporal scales of
the boundary controls to smooth out the oscillations.

\noindent$\bullet$\,\, At the final step, the speed $c$ is found
by (\ref{BASIC}) with the help of numerical differentiation by the
central finite difference formula.

\subsubsection*{Numerical results} \label{sec:test}
\begin{figure}[ht]
\centering
\includegraphics[width=0.75\textwidth]{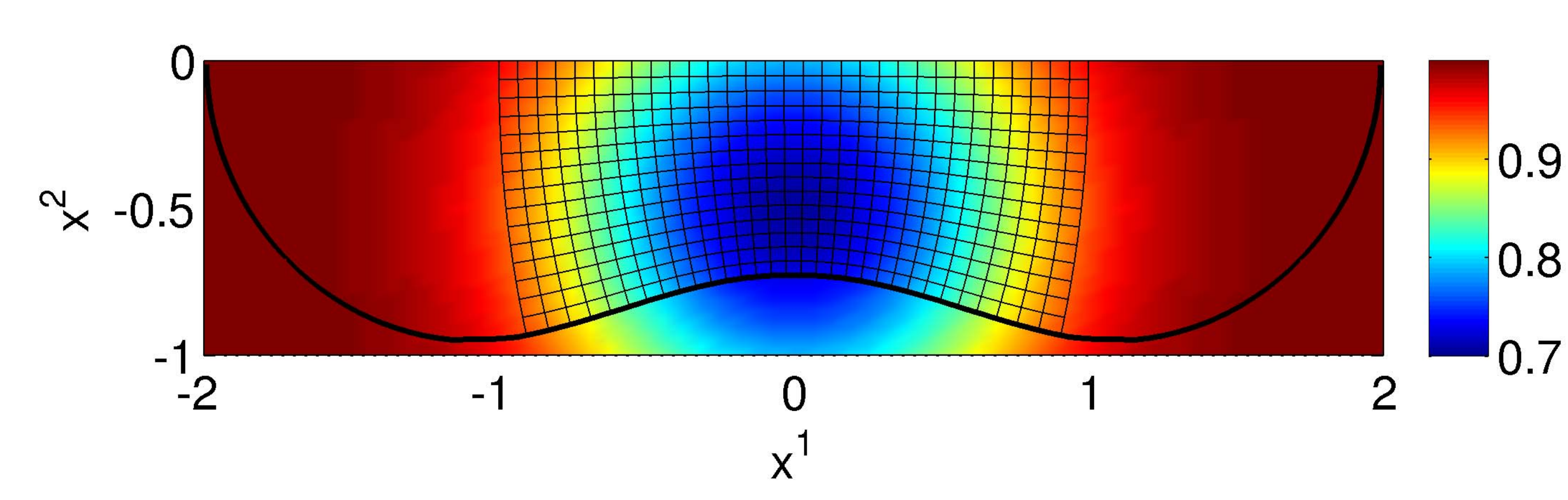}
\caption[]{Test 1. Speed of sound $c(x)$ (color), exact semigeodesic coordinates (mesh), and wave front at $t = T$ (line).}
\label{fig:case1:domain}
\end{figure}

\noindent{\bf Test 1.}\,\,\, Let
\begin{align}
\label{density_1}
\rho(x^1, x^2) = 1 + a\,g_1(x^1)\,g_2(x^2), \quad g_k(x^k) =
\exp{\left[-\frac{\left(x^k - \bar x^k\right)^2}{2\Delta_k^2}\right]},
\end{align}
where $a = 1$, $\bar x^1 = 0$, $\bar x^2 = -0.5$, $\Delta_1 =
0.5$, $\Delta_2 = 0.5$. The sound speed $c = \rho^{-\frac{1}{2}}$
is shown on Figure \ref{fig:case1:domain} together with exact
semigeodesic coordinates and wave front at $t = T = 1$. We test
the recovering procedure for two rather different spatial bases
(the temporal basis (\ref{temporal}) consisting of 16 functions is the same in both cases).

\begin{figure}[ht]
\includegraphics[width=0.47\textwidth]{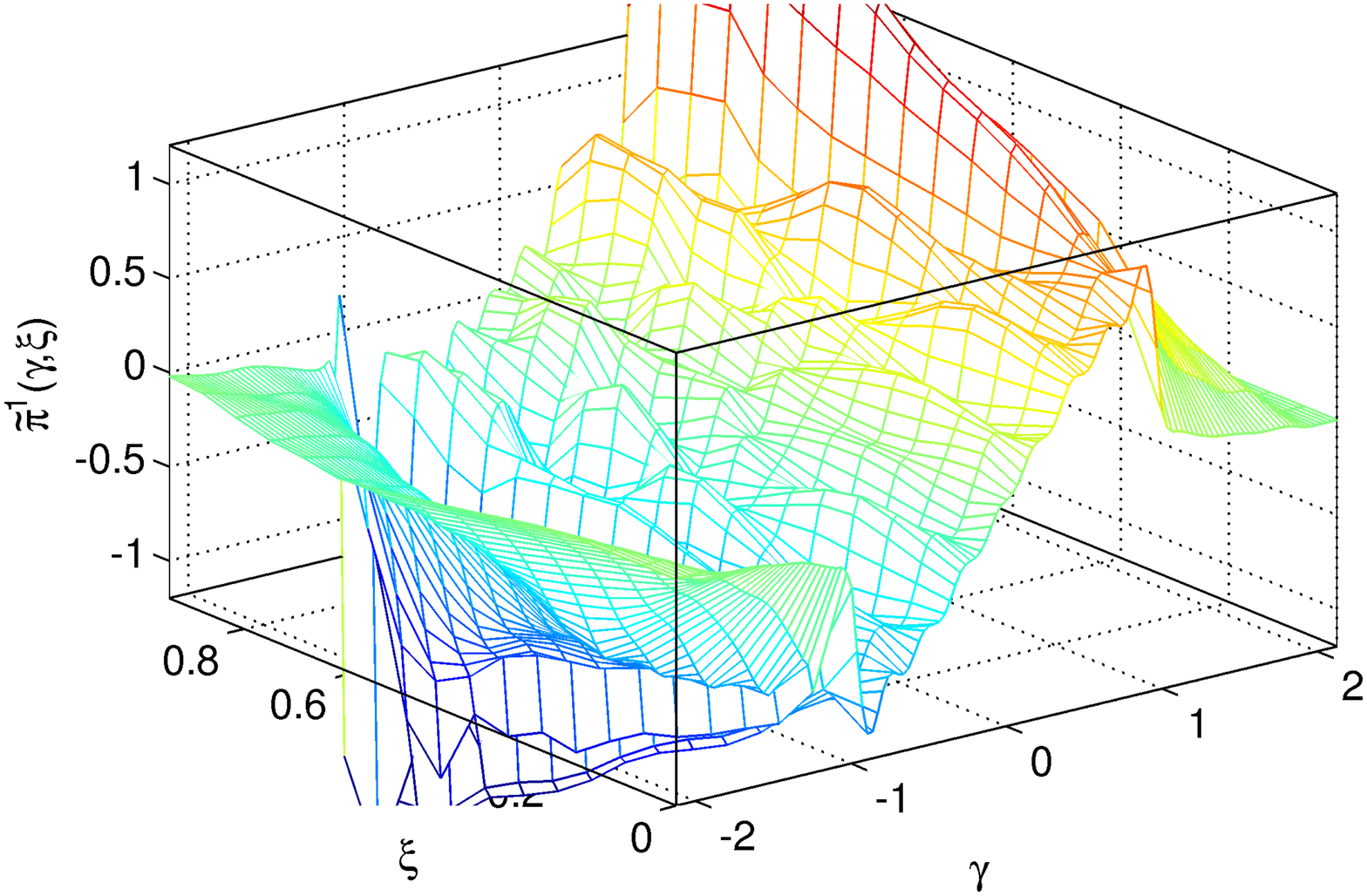}\hfill
\includegraphics[width=0.47\textwidth]{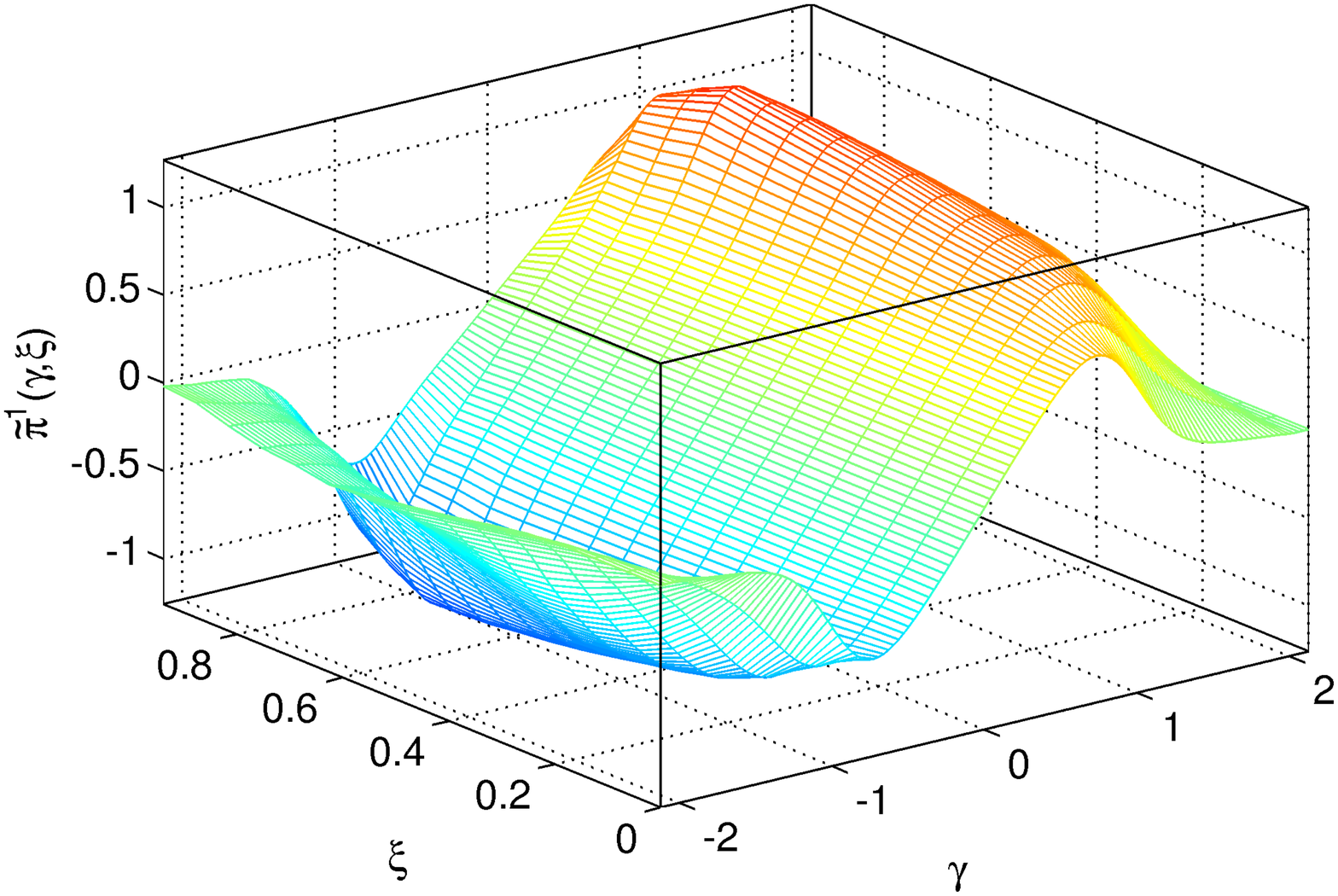}
\caption[]{Test 1. Image $\tilde \pi^1(\gamma, \xi)$: expression (\ref{key rel}) (left) and its smoothed version (\ref{convol}) (right).}
\label{fig:case1:image}
\end{figure}

\noindent{\bf a)}\,\,In this subcase, we use spatial basis
composed from localized trigonometric functions (\ref{spatial}). A
typical image of harmonic function $x^1$ is shown in Figure
\ref{fig:case1:image}, where we observe the Gibbs oscillations on
the left plot and the smoothing effect of convolution
(\ref{convol}) on the right one. The condition number of matrix
(\ref{Gramm G1}) for $\xi = T$ is $1.5{\cdot}10^5$ and parameter
of Tikhonov regularization for all linear systems is
$1{\cdot}10^{-5}$. The standard deviations of Gaussian kernels in
(\ref{convol}) for $(\gamma, t)$ are $\sigma_\gamma = 0.1875$ and
$\sigma_t = 0$.

\begin{figure}[ht]
\includegraphics[width=0.47\textwidth]{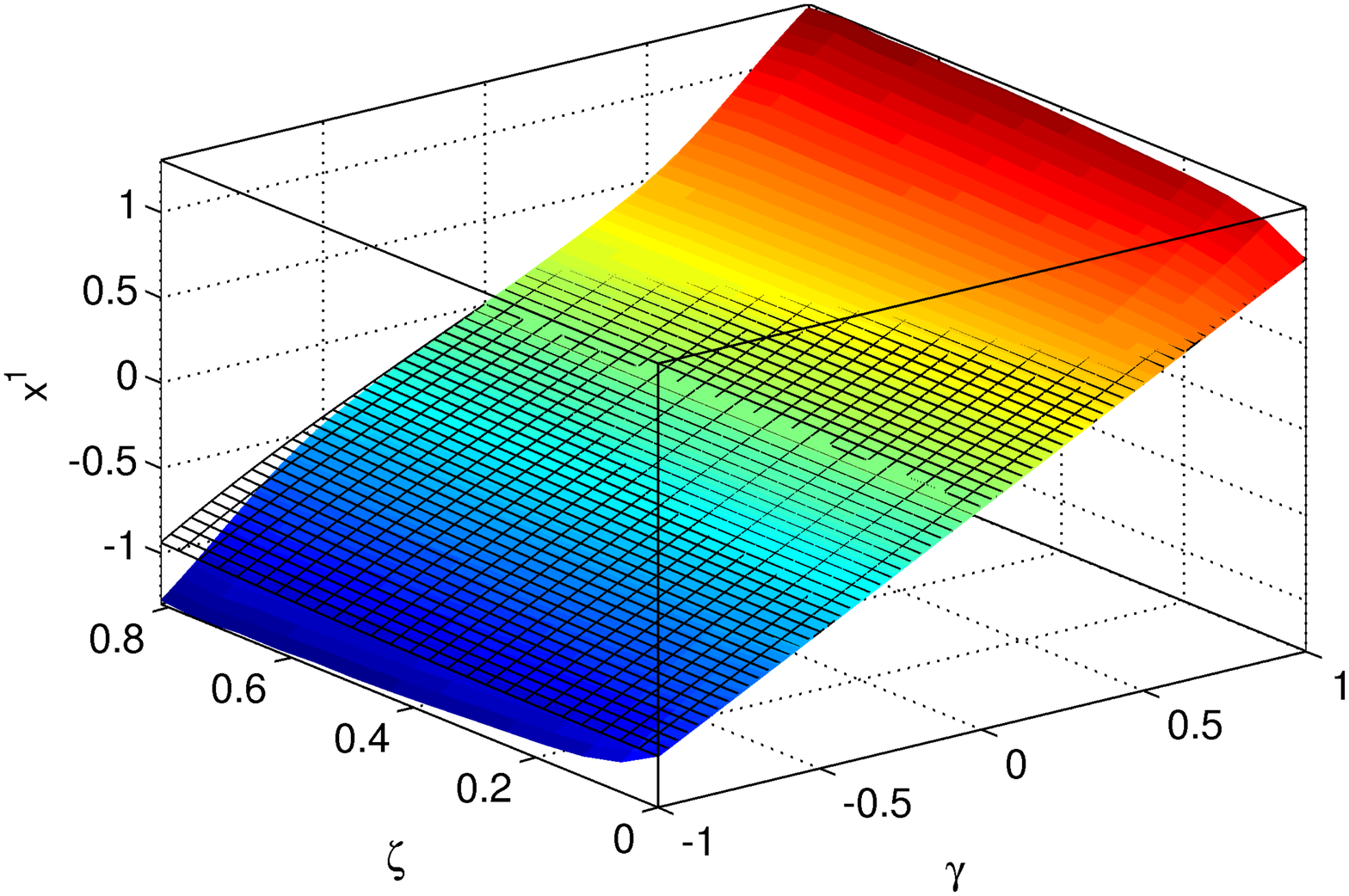}\hfill
\includegraphics[width=0.47\textwidth]{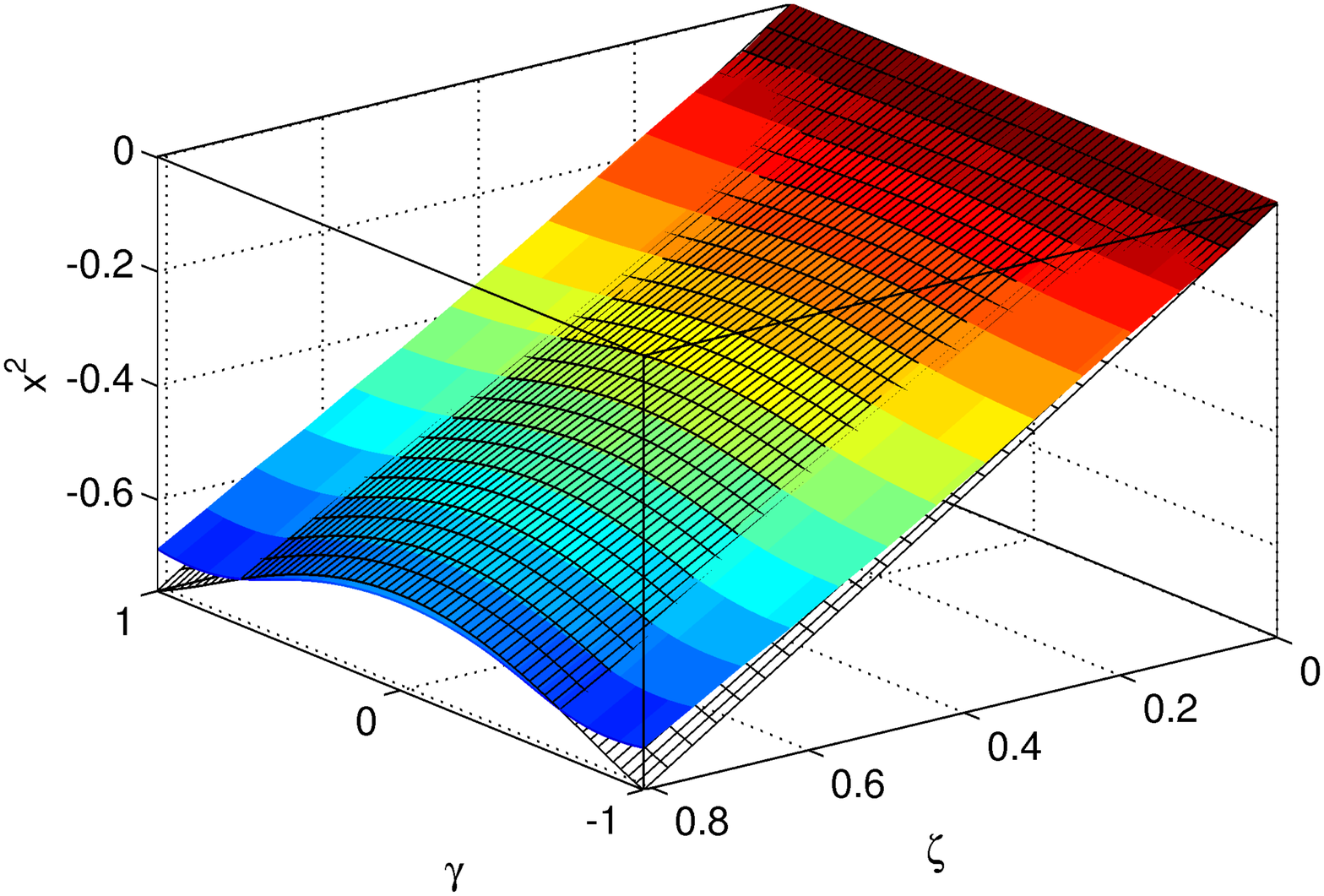}
\caption[]{Test 1. Reconstructed mapping $x = x(\gamma, \xi)$, the exact values are shown by black mesh.}
\label{fig:case1:mapping}
\end{figure}

The mapping $x(\gamma,\xi)$ is shown in Figure
\ref{fig:case1:mapping}. The reconstruction error grows towards
the ends of the localization interval $\gamma\in(-1,1)$ and for
large values of $\xi\approx T$.

\begin{figure}[ht]
\begin{center}
\includegraphics[width=0.325\textwidth]{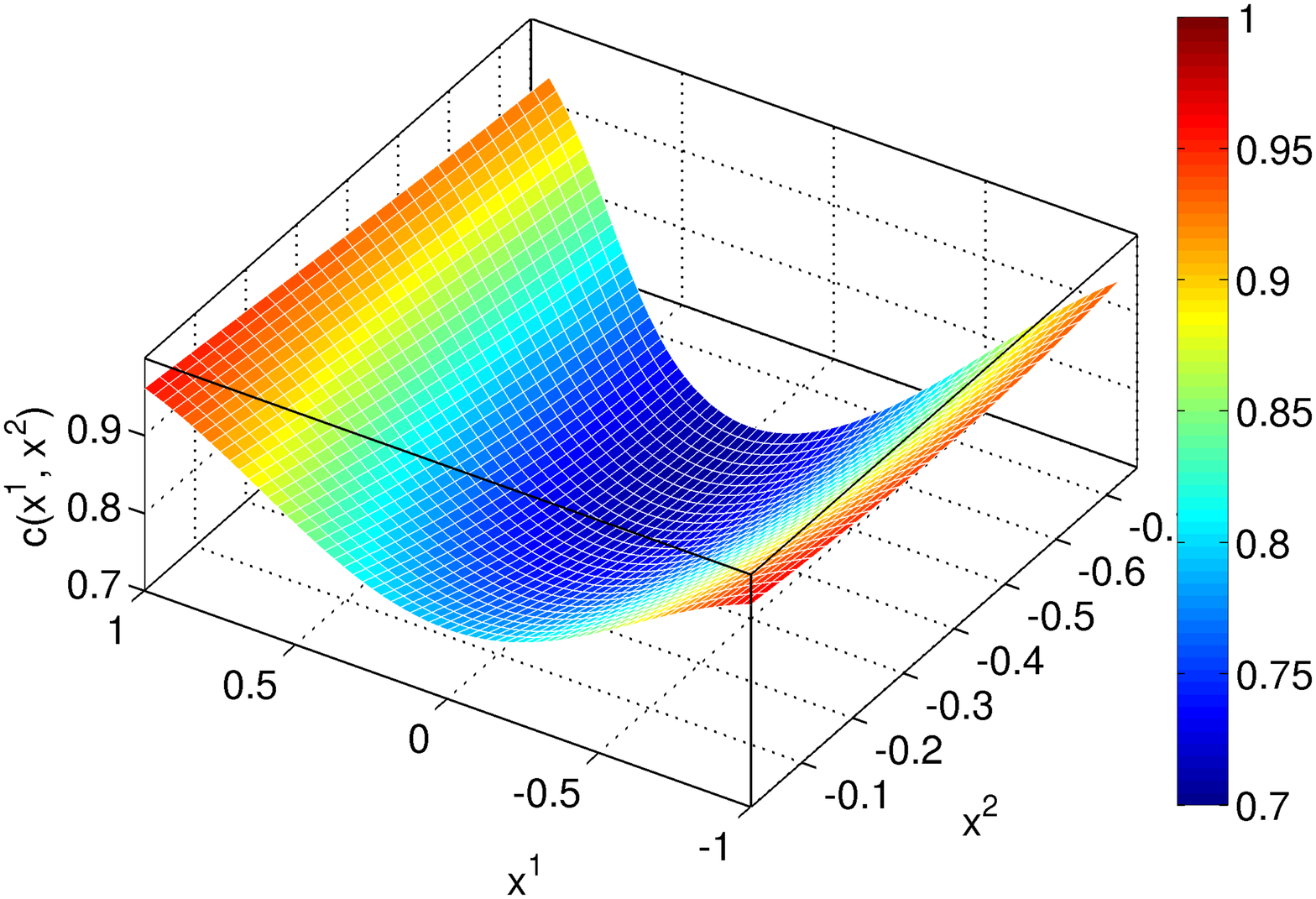}
\includegraphics[width=0.325\textwidth]{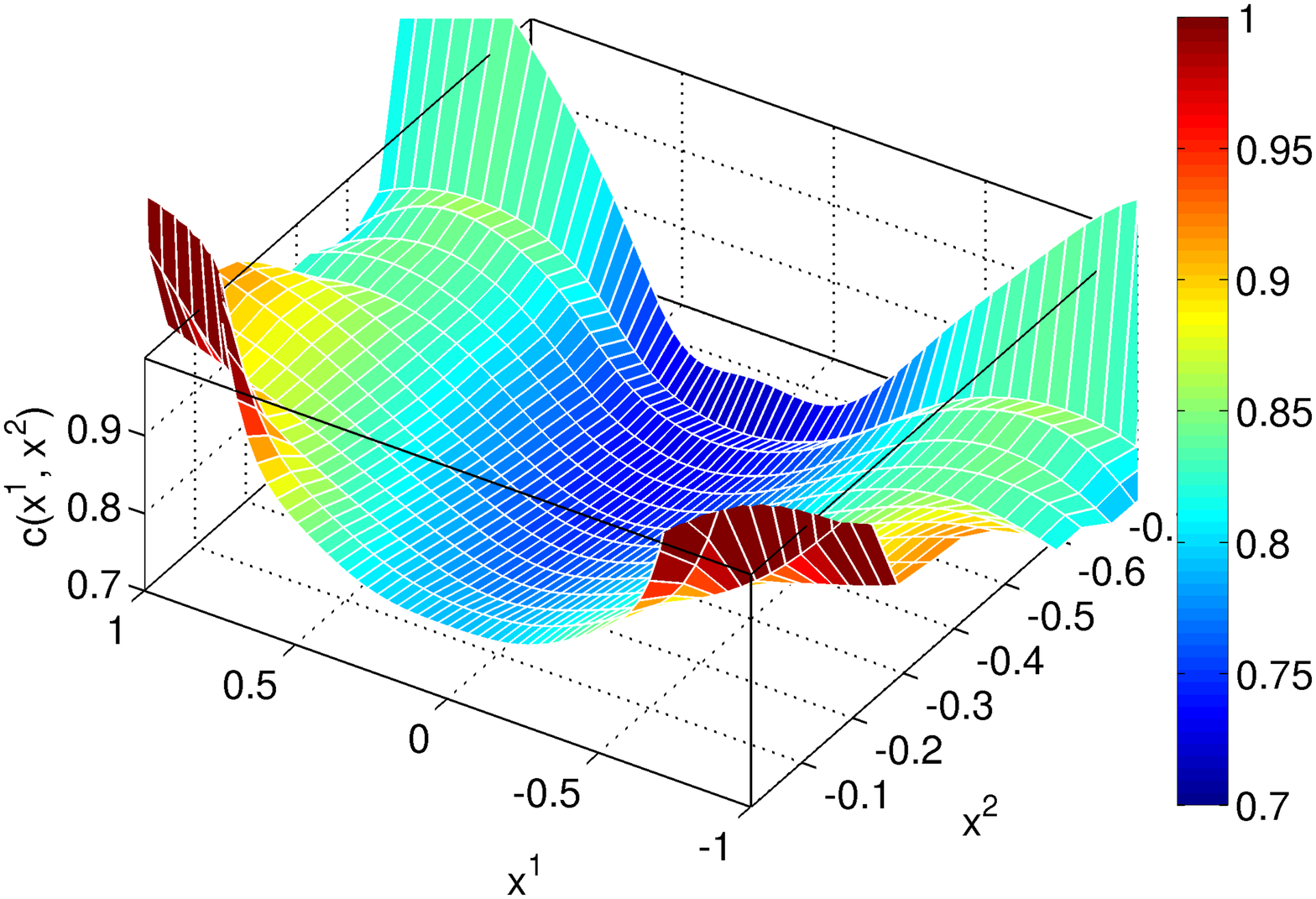}
\includegraphics[width=0.325\textwidth]{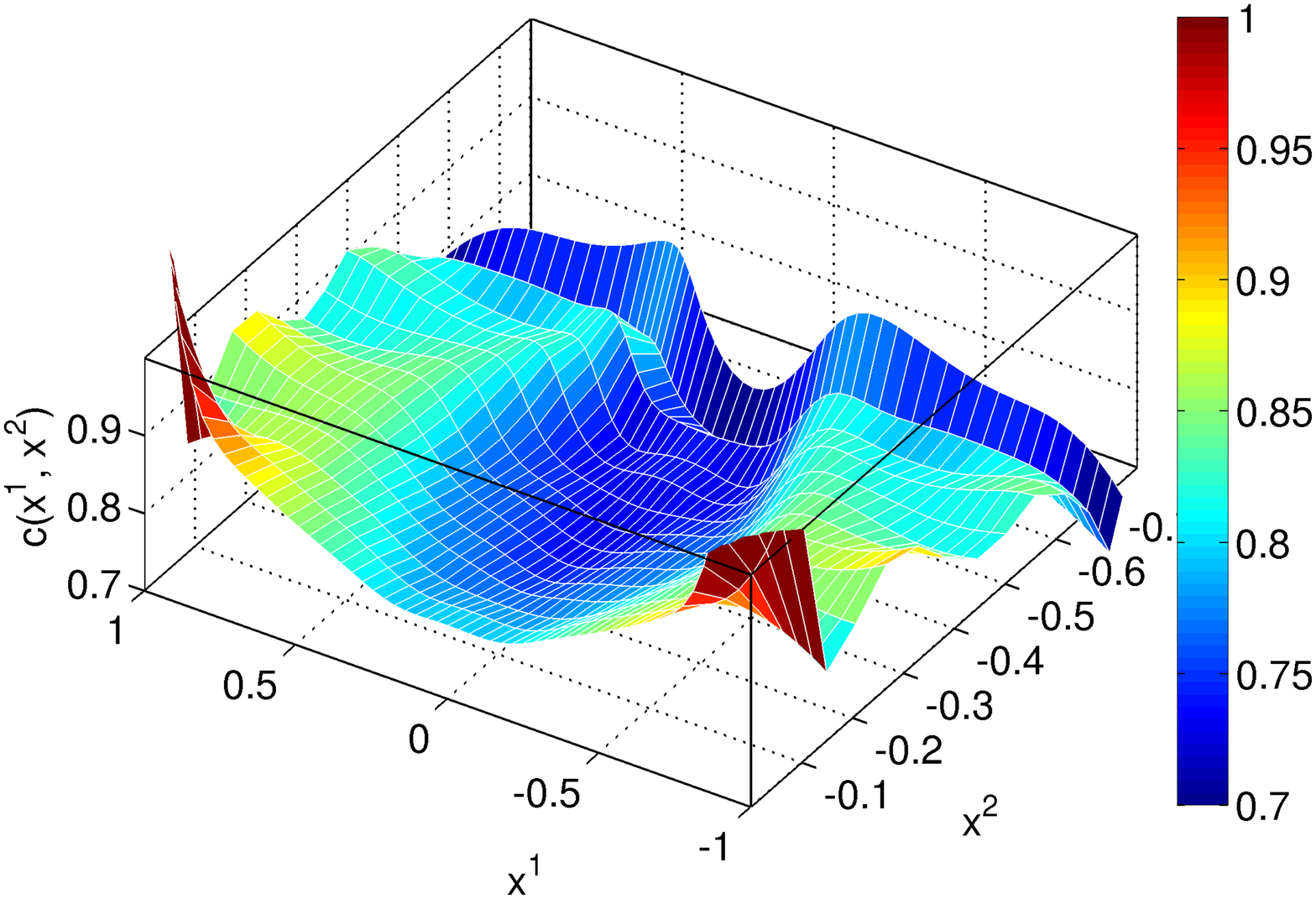}
\caption[]{Test 1. Speed of sound $c(x)$ in the domain filled by waves initiated from $\sigma$:
left - exact values, central - recovered values with trigonometric spatial basis, right - recovered values with tent-like spatial basis.}
\label{fig:case1:vel}
\end{center}
\end{figure}

\begin{figure}[ht]
\begin{center}
\includegraphics[width=0.495\textwidth]{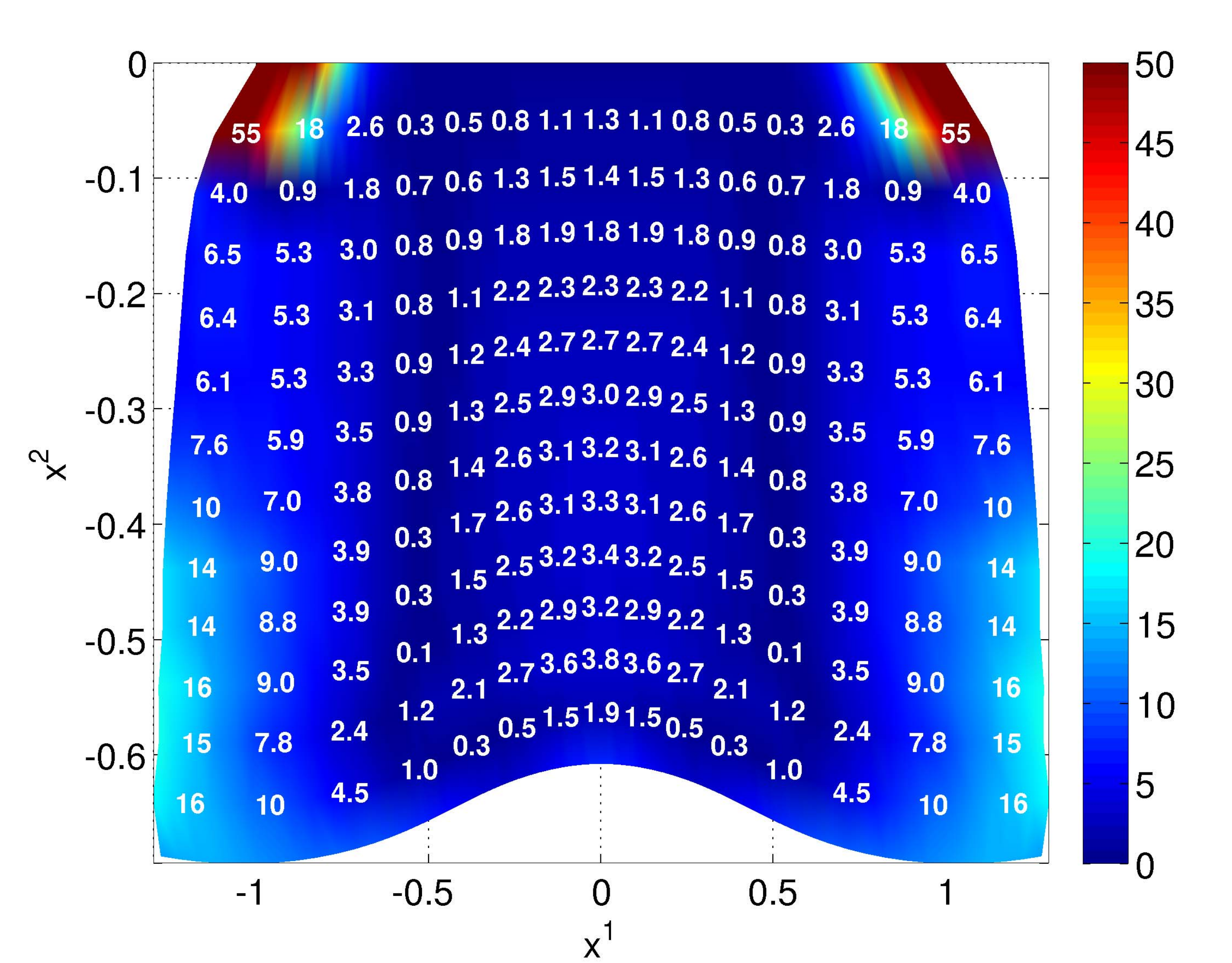}
\includegraphics[width=0.495\textwidth]{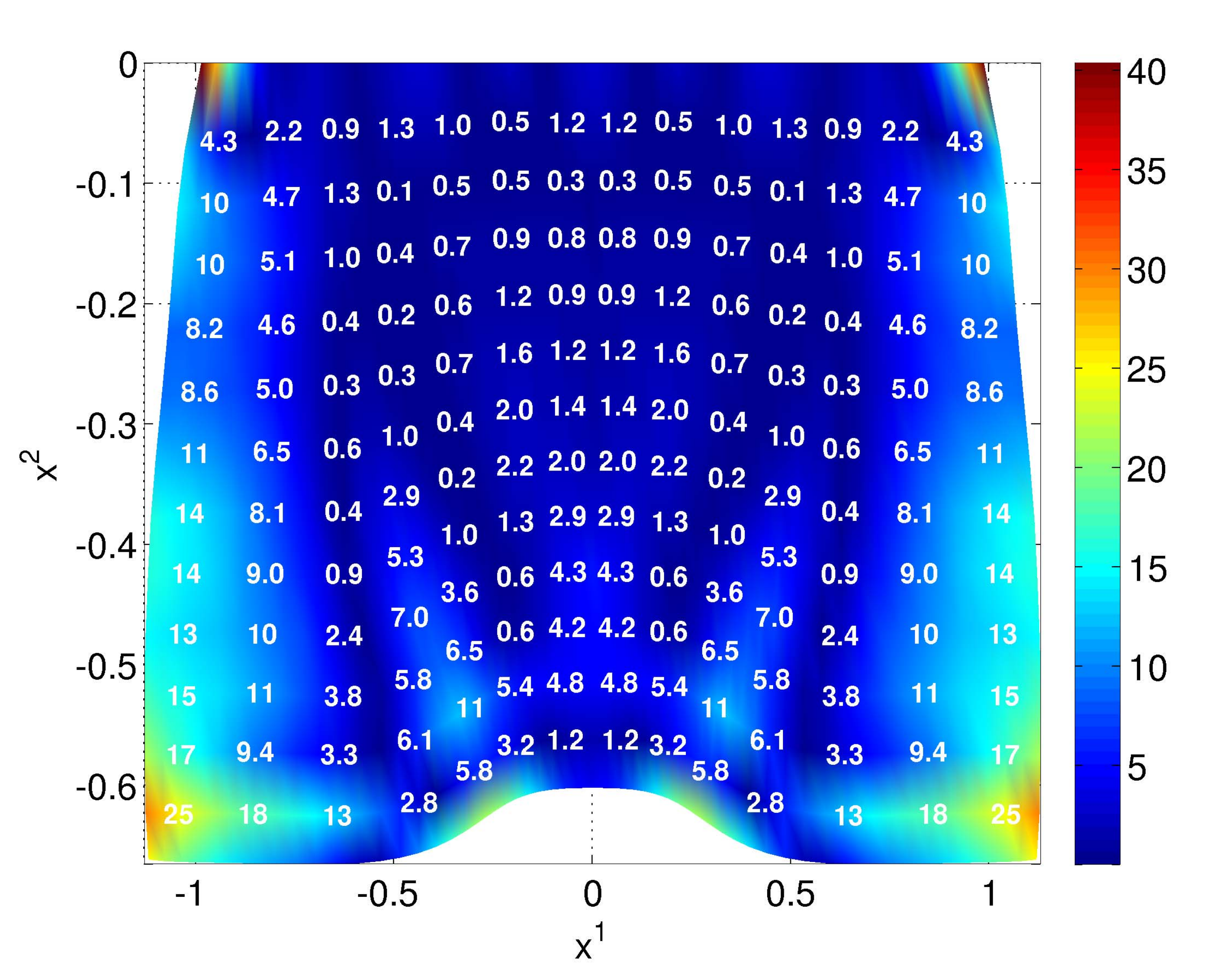}
\caption[]{Test 1. Map of relative errors (in percents) of the recovered sound speed $c(x)$:
left - with trigonometric spatial basis, right - with tent-like spatial basis.}
\label{fig:case1:relerr}
\end{center}
\end{figure}

The end result of the BCM is the sound speed recovered in the
Cartesian coordinates. It is shown in central panel of Figure
\ref{fig:case1:vel}. Relative errors of reconstruction in percents
are shown in left panel of Figure \ref{fig:case1:relerr}. As is
seen, although the reconstruction error quickly grows towards the
ends of the localization interval and for large values of
$\xi\approx T$, in the most part of the domain covered by the
direct rays from the boundary, the relative error does not exceed
a few percents.

\noindent{\bf b)}\,\, Here we use a spatial basis composed from
smooth tent-like functions as in (\ref{temporal}). The condition
number of matrix (\ref{Gramm G1}) for $\xi = T$ is
$5.9{\cdot}10^3$ and parameter of Tikhonov regularization for all
linear systems is $1{\cdot}10^{-6}$. The standard deviations of
Gaussian kernels in (\ref{convol}) for $(\gamma, t)$ are
$\sigma_\gamma = 0.125$ and $\sigma_t = 0.0625$. The recovered
speed of sound is shown in right panel of Figure
\ref{fig:case1:vel} while its relative errors are shown in right
panel of Figure \ref{fig:case1:relerr}.

We may conclude that both of these bases provide similar quality
of reconstruction of the order of several percents in most part of
the domain. The advantages of the tent-like basis are smaller
condition number of the system matrix and the same spatial scale
of all basis functions. The effect of lower accuracy of
reconstruction along the lateral boundaries in the case (b) is due
to narrower support (smaller value of $L$) of the boundary
controls compared to the case (a).

\medskip

\noindent{\bf Test 2.}\,\,\, For the second test, we take
\begin{align}
\rho(x^1, x^2) = 1 - 0.5 x^2 + 0.0625 \left(x^1\right)^2 -
a\,g_1(x^1)\,\frac{\partial g_2(x^2)}{\partial x^2},\nonumber
\end{align}
where $a = 0.25$, $\bar x^1 = 0$, $\bar x^2 = -0.5$, $\Delta_1 =
0.5$, $\Delta_2 = 0.25$. The corresponding sound speed has a
background value $1$ and two variations of the order 30\% of its
boundary value.

\begin{figure}[ht]
\centering
\includegraphics[width=0.54\textwidth]{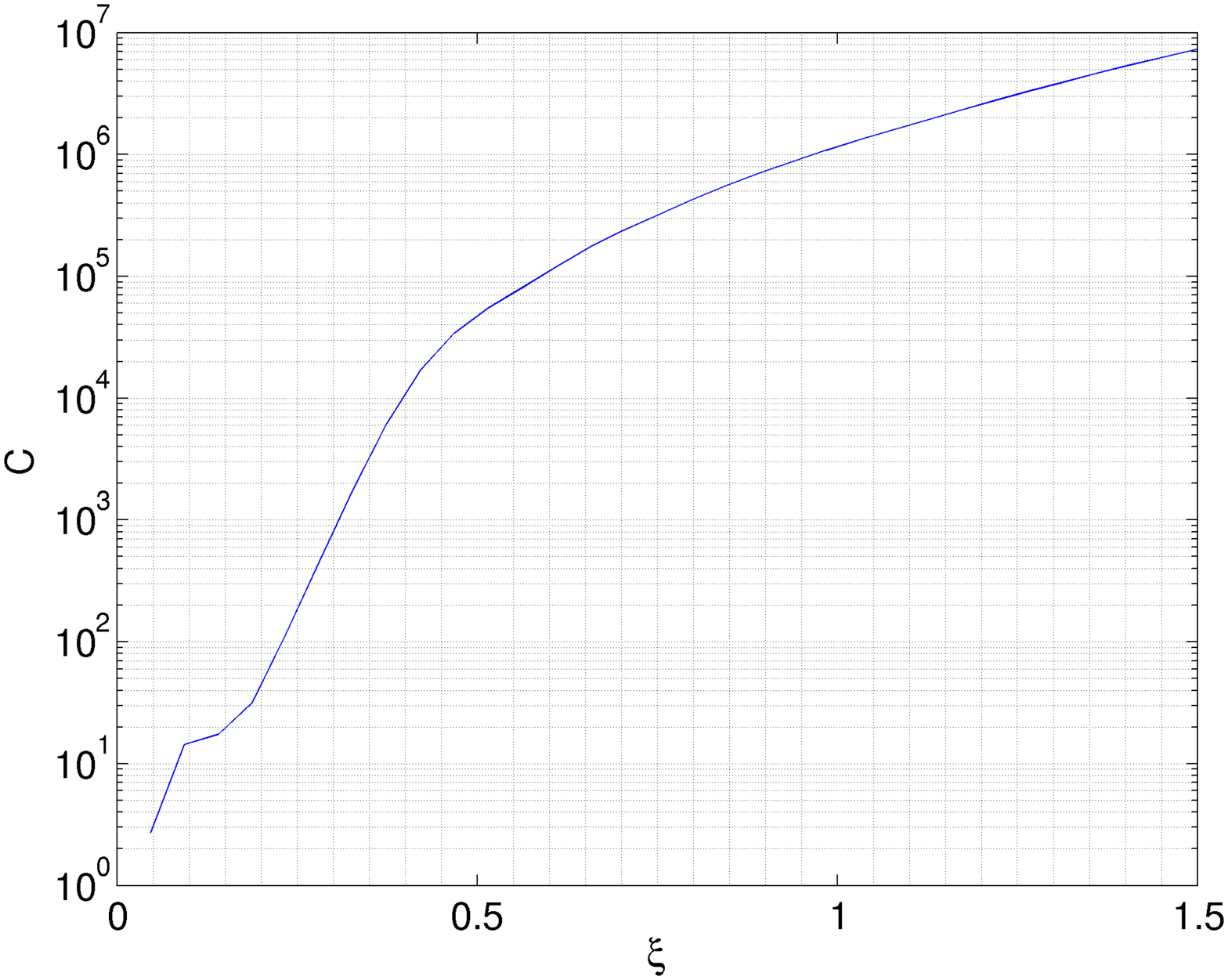}
\caption[]{Test 2. The condition number of matrix (\ref{Gramm G1}) of scalar products as function of the probing time $\xi$.}
\label{fig:case2:cond}
\end{figure}

We use $T = 1.5$ and the basis with 16 spatial (trigonometric) and 32 temporal
functions. The condition number of matrix (\ref{Gramm G1}) is
shown in Figure \ref{fig:case2:cond}; it grows as $\xi^4$. This is
a consequence of the strong ill-posedness of the inverse problem
under consideration, and such a growth constrains the maximal
depth of reconstruction (determined by errors in right hand sides
(\ref{B})), which is possible for the given part of the boundary.
In computations, the parameter of Tikhonov regularization for all
linear systems is fixed and equal to $1{\cdot}10^{-4}$. The
standard deviations of Gaussian kernels in (\ref{convol}) for
$(\gamma, t)$ are $\sigma_\gamma = 0.1875$ and $\sigma_t =
4.6875{\cdot}10^{-2}$. For $\xi\approx T$, the error in the
expansion coefficients $a^\xi_\alpha$ grows up and we had to
increase $\sigma_\gamma$ to the value $0.5$ for smoothing out the
large scale fake oscillations from low spatial harmonics. Such an
over-smoothing reduces the accuracy of the recovering for
$\xi\approx T$.

\begin{figure}[ht]
\begin{center}
\includegraphics[width=0.495\textwidth]{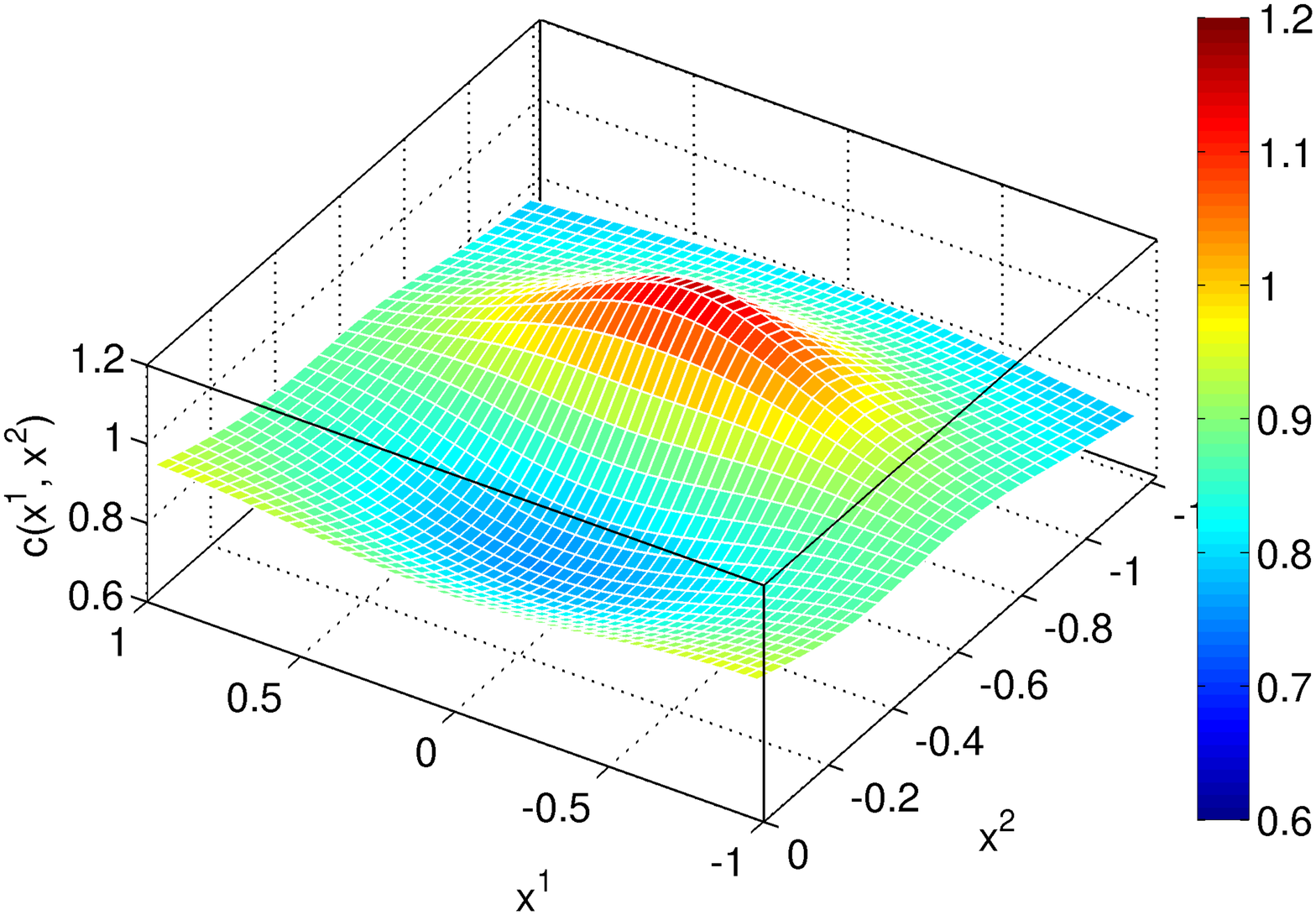}
\includegraphics[width=0.495\textwidth]{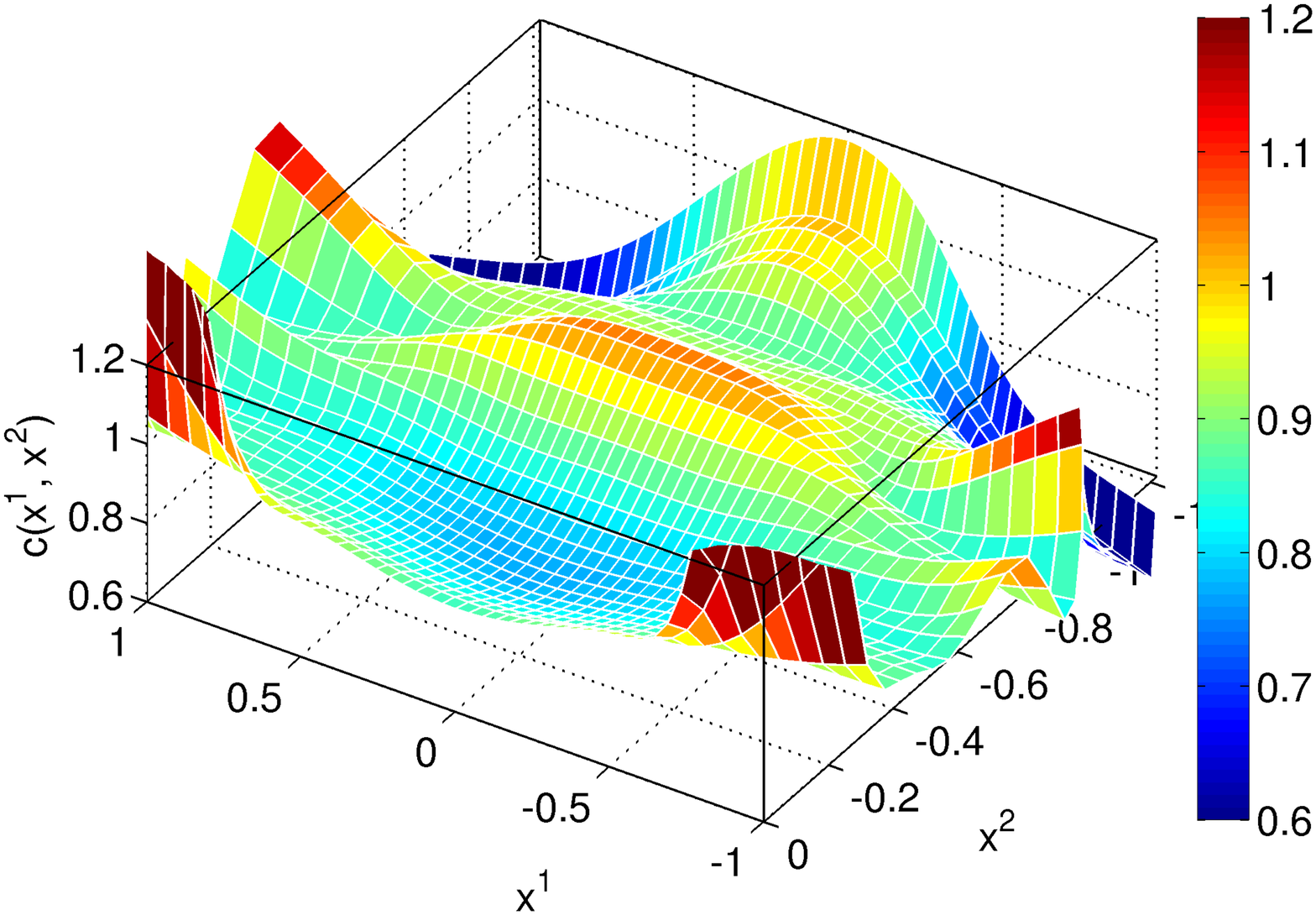}
\caption[]{Test 2. Speed of sound $c(x)$ in the domain filled by waves initiated from $\sigma$:
left - exact values, right - recovered values.}
\label{fig:case2:vel}
\end{center}
\end{figure}

\begin{figure}[ht]
\begin{center}
\includegraphics[width=0.58\textwidth]{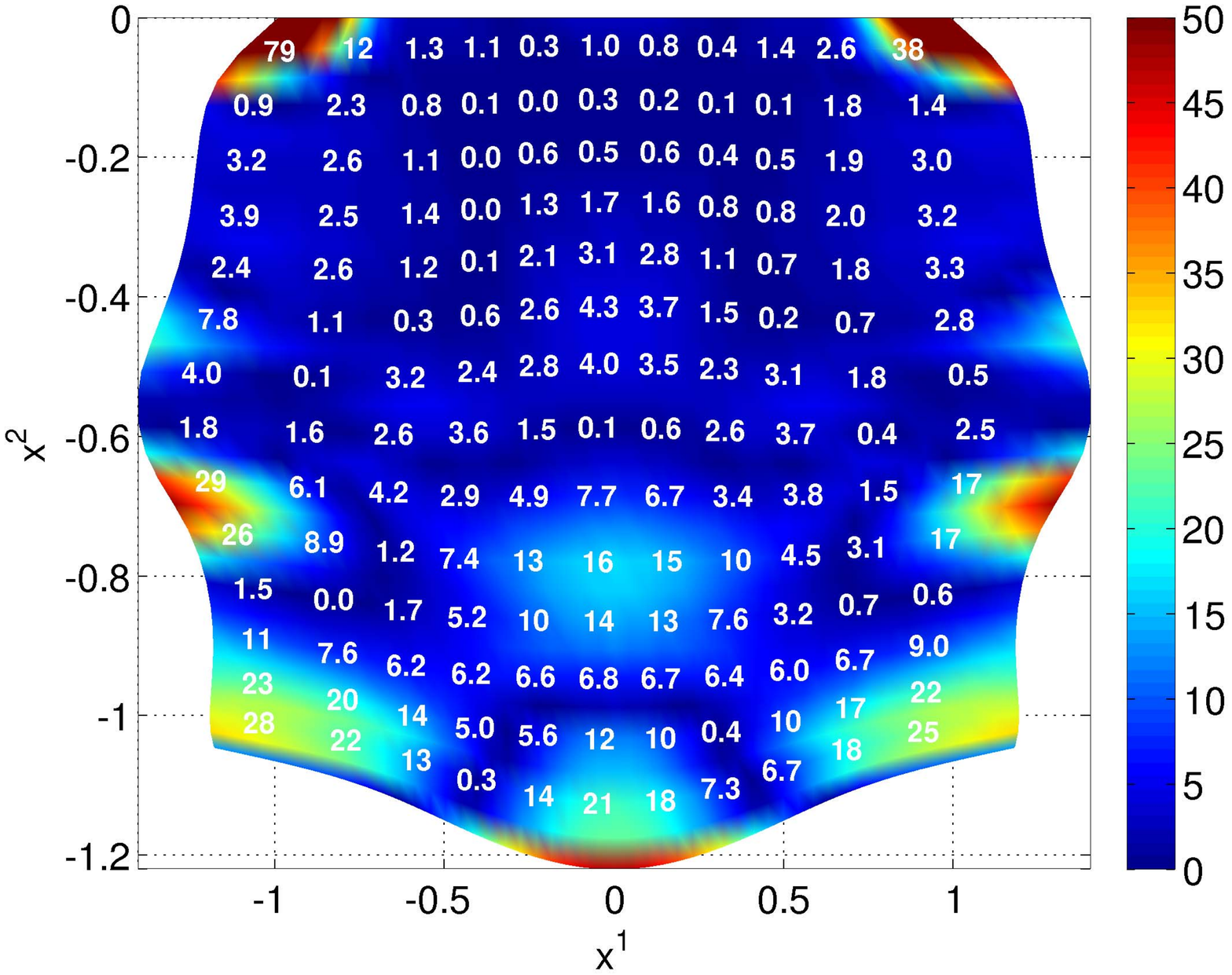}
\caption[]{Test 2. Map of relative errors of the recovered sound speed in percents.}
\label{fig:case2:relerr}
\end{center}
\end{figure}

The recovered speed of sound in the Cartesian coordinates is shown
in Figure \ref{fig:case2:vel}, and relative errors of
reconstruction in percents are shown in Figure
\ref{fig:case2:relerr}. Thus, in the most part of the domain
$\Omega_\sigma$ covered by the direct rays coming from $\sigma$,
the relative error does not exceed a few percents.
\medskip

\noindent{\bf Test 3.}\,\,\, Here we take
\begin{align}
\rho(x^1, x^2) = 1 - 0.5 x^2 + 0.0625 \left(x^1\right)^2 + a\,g_1(x^1)\,\left(1 - x^2\right)\,
\frac{\partial g_2(x^2)}{\partial x^2},\nonumber
\end{align}
where $a = 0.25$, $\bar x^1 = 0$, $\bar x^2 = -0.5$, $\Delta_1 =
0.5$, $\Delta_2 = 0.25$. In contrast to the case 2, the
corresponding speed of sound has rather strong variations (the
ratio of maximum to minimum value is about $2.5$).

Again, we take $T = 1.5$ and use the basis with 16 spatial and 32
temporal functions. The condition number of matrix (\ref{Gramm
G1}) for $\xi = T$ is $2{\cdot}10^5$ and in calculations the
parameter of Tikhonov regularization for all linear systems is
$1{\cdot}10^{-4}$. The standard deviations of Gaussian kernels in
(\ref{convol}) for $(\gamma, t)$ are $\sigma_\gamma = 0.1875$ and
$\sigma_t = 4.6875{\cdot}10^{-2}$. Again, for $\xi\approx T$ we
had to increase $\sigma_\gamma$ to value $0.625$ to smooth out
large scale fake oscillations from the low spatial harmonics.

\begin{figure}[ht]
\begin{center}
\includegraphics[width=0.495\textwidth]{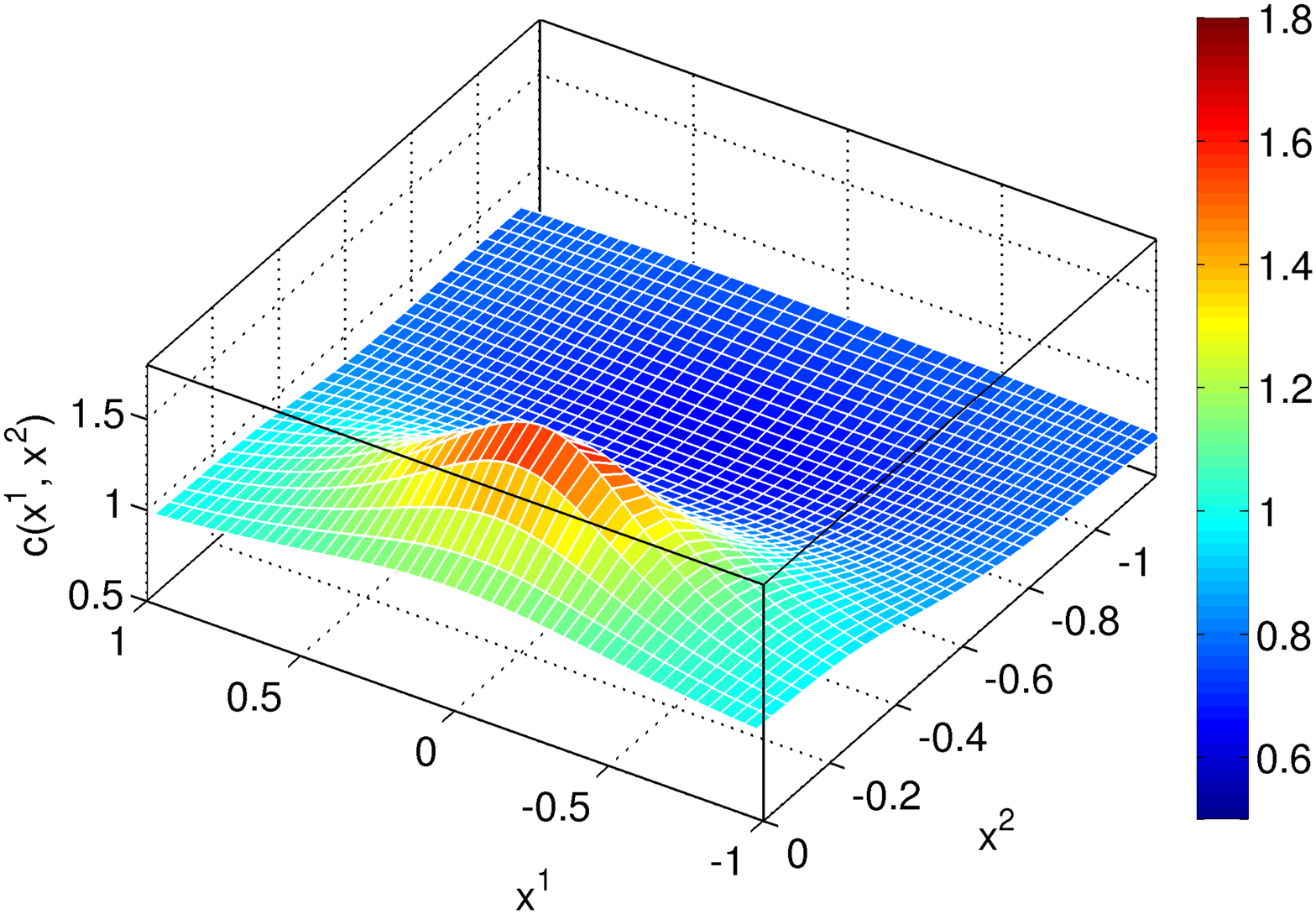}
\includegraphics[width=0.495\textwidth]{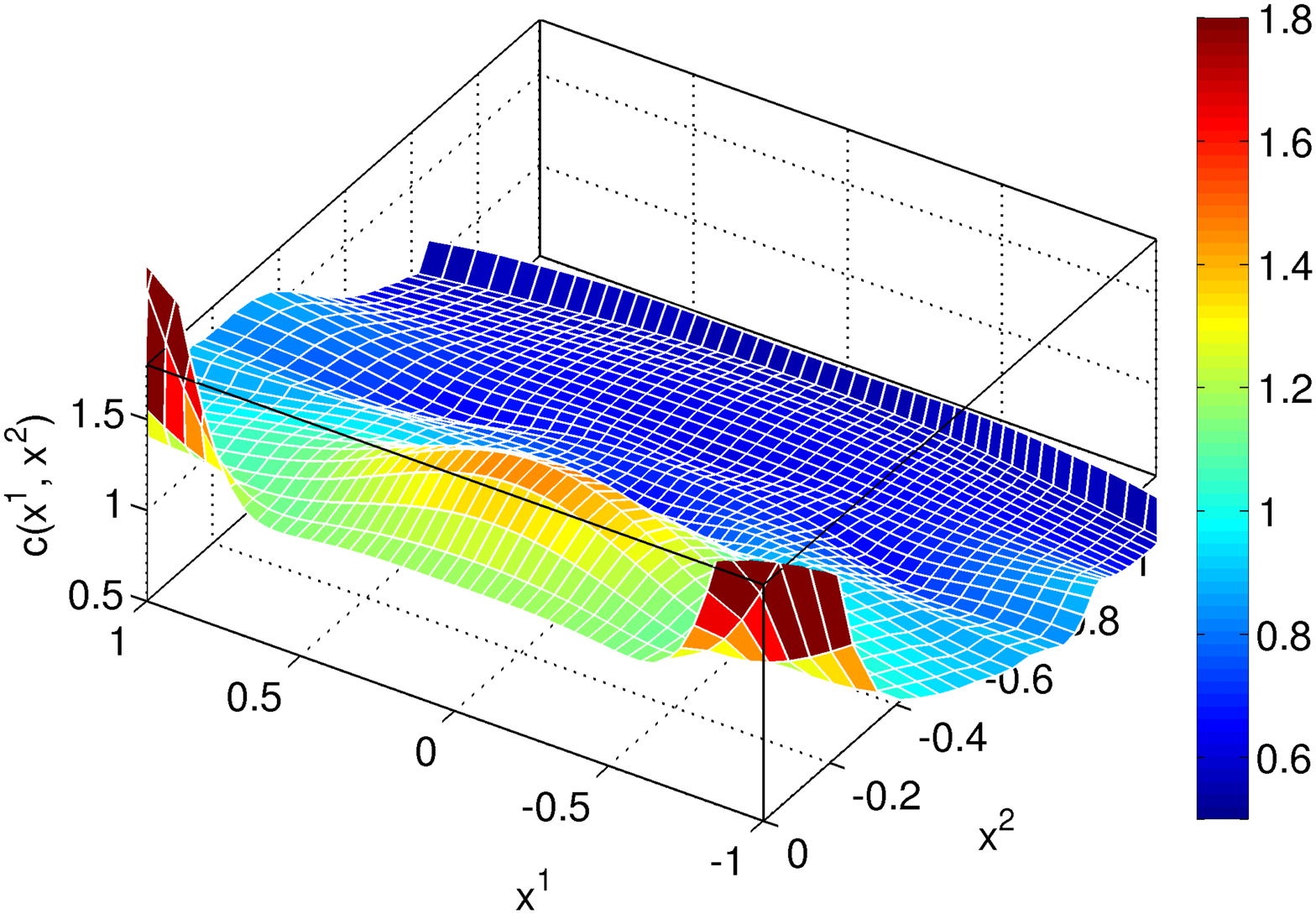}
\caption[]{Test 3. Speed of sound $c(x)$ in the domain filled by waves initiated from $\sigma$:
left - exact values, right - recovered values.}
\label{fig:case3:vel}
\end{center}
\end{figure}

\begin{figure}[ht]
\begin{center}
\includegraphics[width=0.57\textwidth]{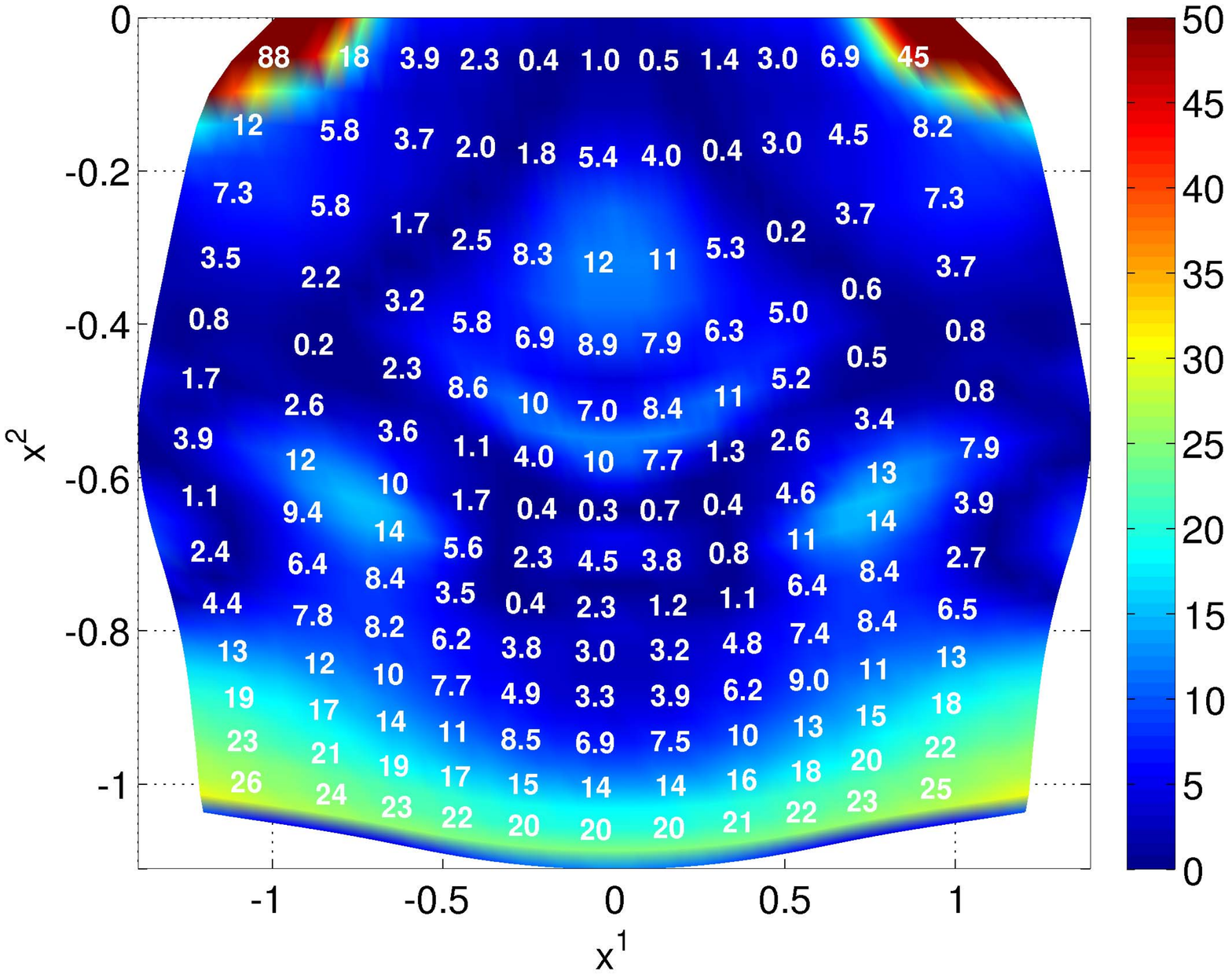}
\caption[]{Test 3. Map of relative errors of the recovered sound speed in percents.}
\label{fig:case3:relerr}
\end{center}
\end{figure}

The recovered speed of sound in the Cartesian coordinates is shown
in Figure \ref{fig:case3:vel}, and relative errors of
reconstruction in percents are shown in Figure \ref{fig:case3:relerr}.

\medskip

\noindent{\bf Test 4.}\,\,\, To test the recovering algorithm for
the case of sound speed quickly varying along the boundary we set
\begin{align}
\label{density_4}
\rho(x^1, x^2) &= 1 - a\,g_2(z^2)\,\frac{\partial g_1(z^1)}{\partial x^1},\nonumber\\
z^1 &=   \cos(\phi) x^1 + \sin(\phi) (x^2 + 0.25),\nonumber\\
z^2 &= - \sin(\phi) x^1 + \cos(\phi) (x^2 + 0.25),\nonumber
\end{align}
where $a = 0.25$, $\bar x^1 = 0$, $\bar x^2 = 0$, $\Delta_1 = 0.375$, $\Delta_2 = 0.25$, $\phi = \pi / 12$.

We take $T = 1$ and use the basis with 16 spatial and 32 temporal
functions. The condition number of matrix (\ref{Gramm G1}) for
$\xi = T$ is $1.28{\cdot}10^5$ and in calculations the parameter
of Tikhonov regularization for all linear systems is
$1{\cdot}10^{-4}$. The standard deviations of Gaussian kernels in
(\ref{convol}) for $(\gamma, t)$ are $\sigma_\gamma = 0.125$ and
$\sigma_t = 3.125{\cdot}10^{-2}$. Again, for $\xi\approx T$ we had
to increase $\sigma_\gamma$ to the value $0.25$ in order to
decrease large scale fake oscillations from the low spatial
harmonics.

\begin{figure}[ht]
\begin{center}
\includegraphics[width=0.325\textwidth]{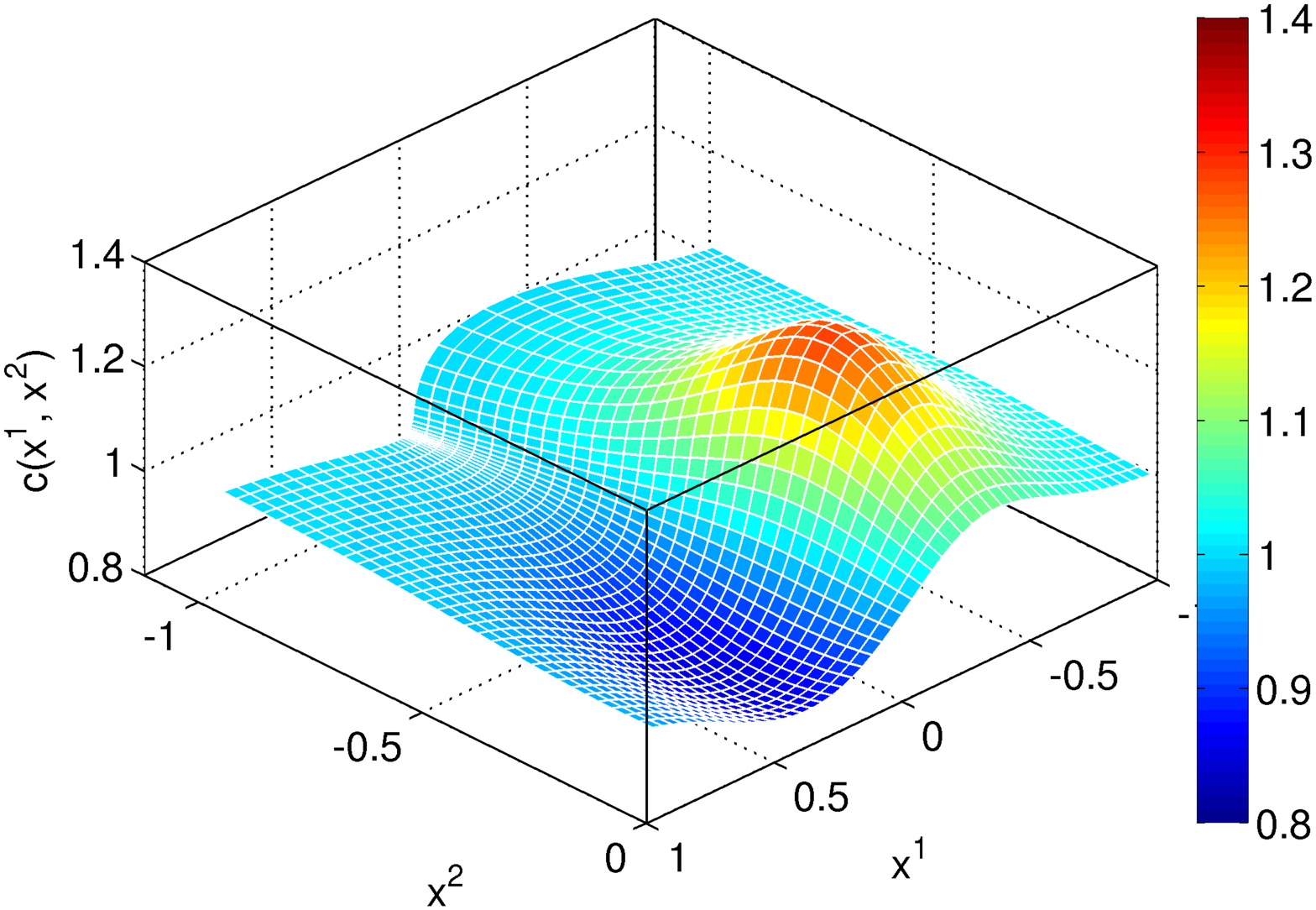}
\includegraphics[width=0.325\textwidth]{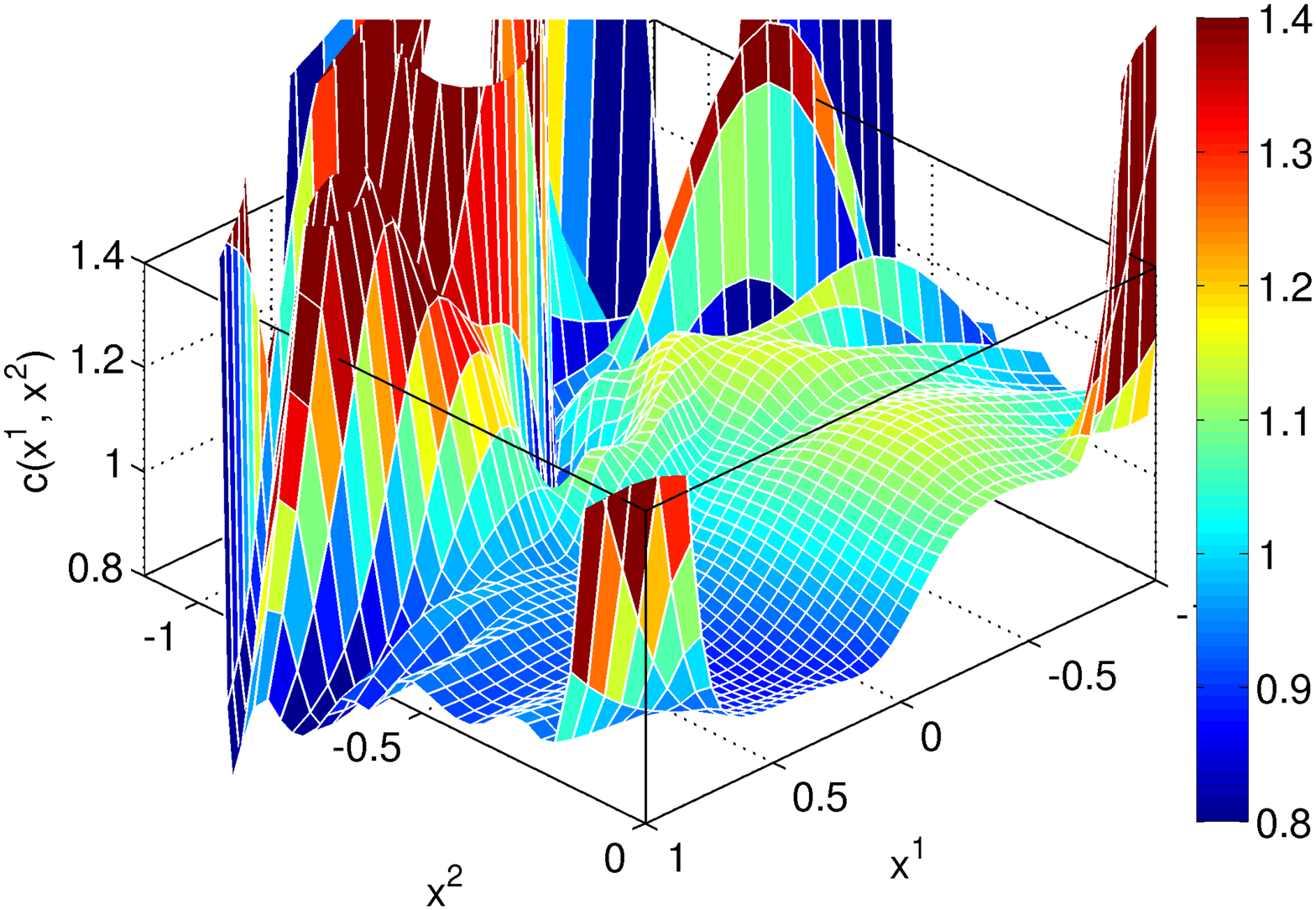}
\includegraphics[width=0.325\textwidth]{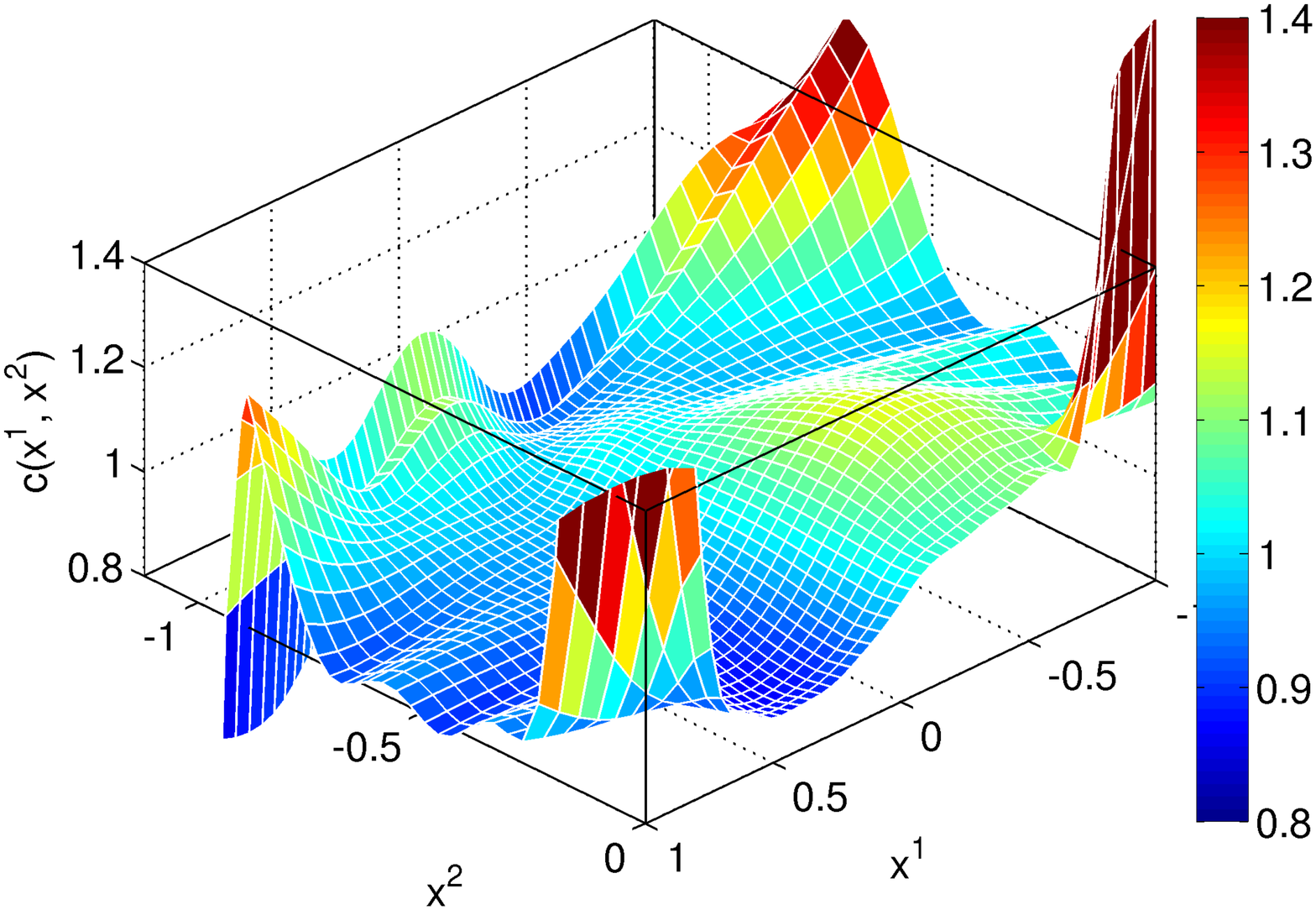}
\caption[]{Test 4. Speed of sound $c(x)$ in the domain filled by waves initiated from $\sigma$:
left - exact values, central - recovered values, right - pseudo-recovered values.}
\label{fig:case4:vel}
\end{center}
\end{figure}

\begin{figure}[ht]
\begin{center}
\includegraphics[width=0.495\textwidth]{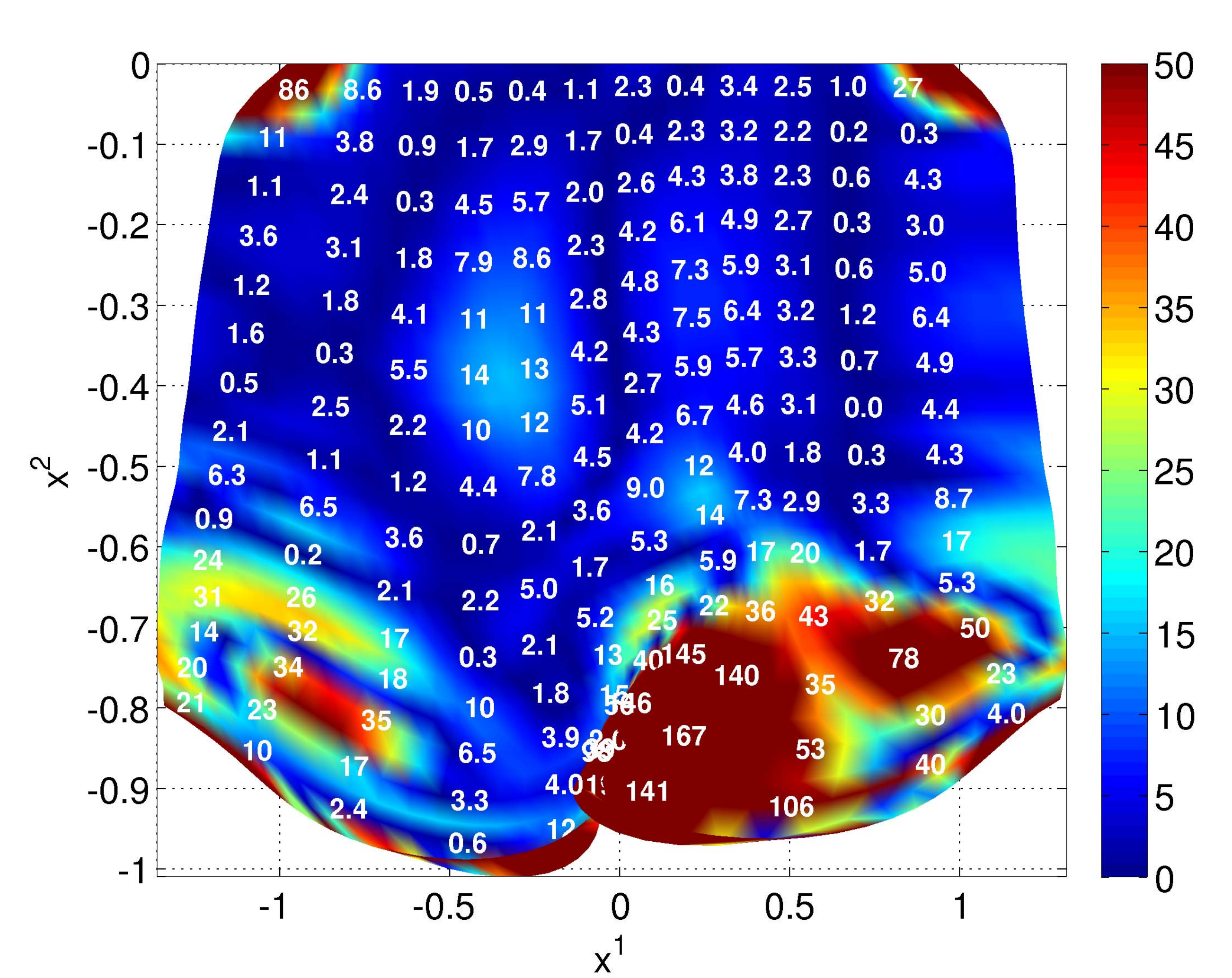}
\includegraphics[width=0.495\textwidth]{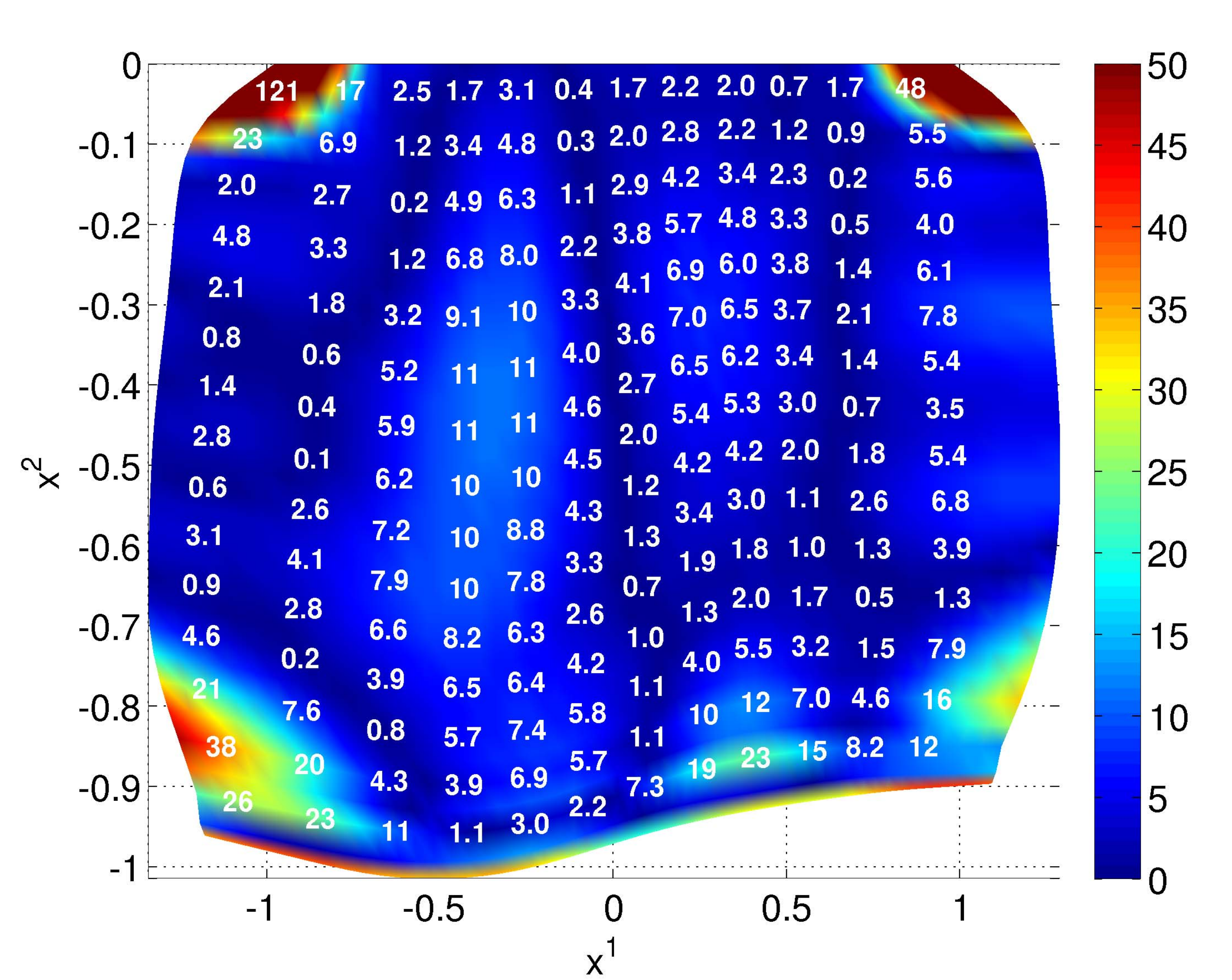}
\caption[]{Test 4. Map of relative errors (in percents) of the recovered sound speed $c(x)$:
left - usual reconstruction, right - pseudo-reconstruction.}
\label{fig:case4:relerr}
\end{center}
\end{figure}

The recovered speed in the Cartesian coordinates is shown in
Figure \ref{fig:case4:vel}, and relative errors of reconstruction
in percents are shown in Figure \ref{fig:case4:relerr}. We observe
rather large oscillations in the speed values emerging at $x^2 = -
0.6$, and to demonstrate the origin of these oscillations we also
show the results of a pseudo-reconstruction. The latter means the
use of a conventional recovering procedure, in which all the
matrix elements (products $(u^{f_i},u^{f_j})$) are computed via
the solutions $u^{f_i}$ found by solving the forward problem with
the given (known) speed profile. Such products are much more
accurate than the ones found via the inverse data, and therefore
the errors in the expansion coefficients in (\ref{key rel}) are
greatly reduced in the pseudo-reconstruction. This leads to much
better quality of recovering far from the boundary and clearly
shows the effect of ill-posedness of the reconstruction procedure.
\medskip

\begin{figure}[ht]
\begin{center}
\includegraphics[width=0.49\textwidth]{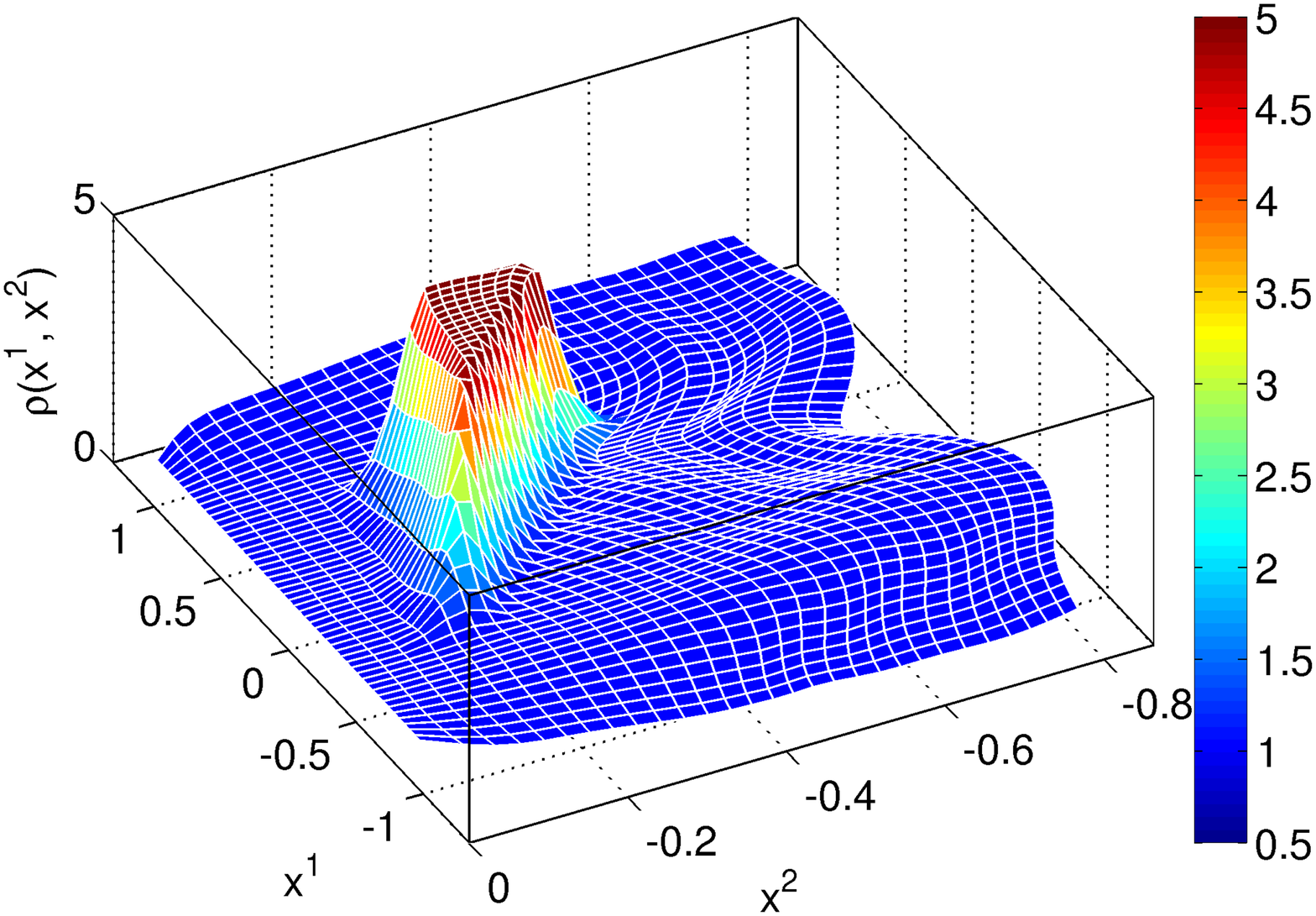}
\includegraphics[width=0.49\textwidth]{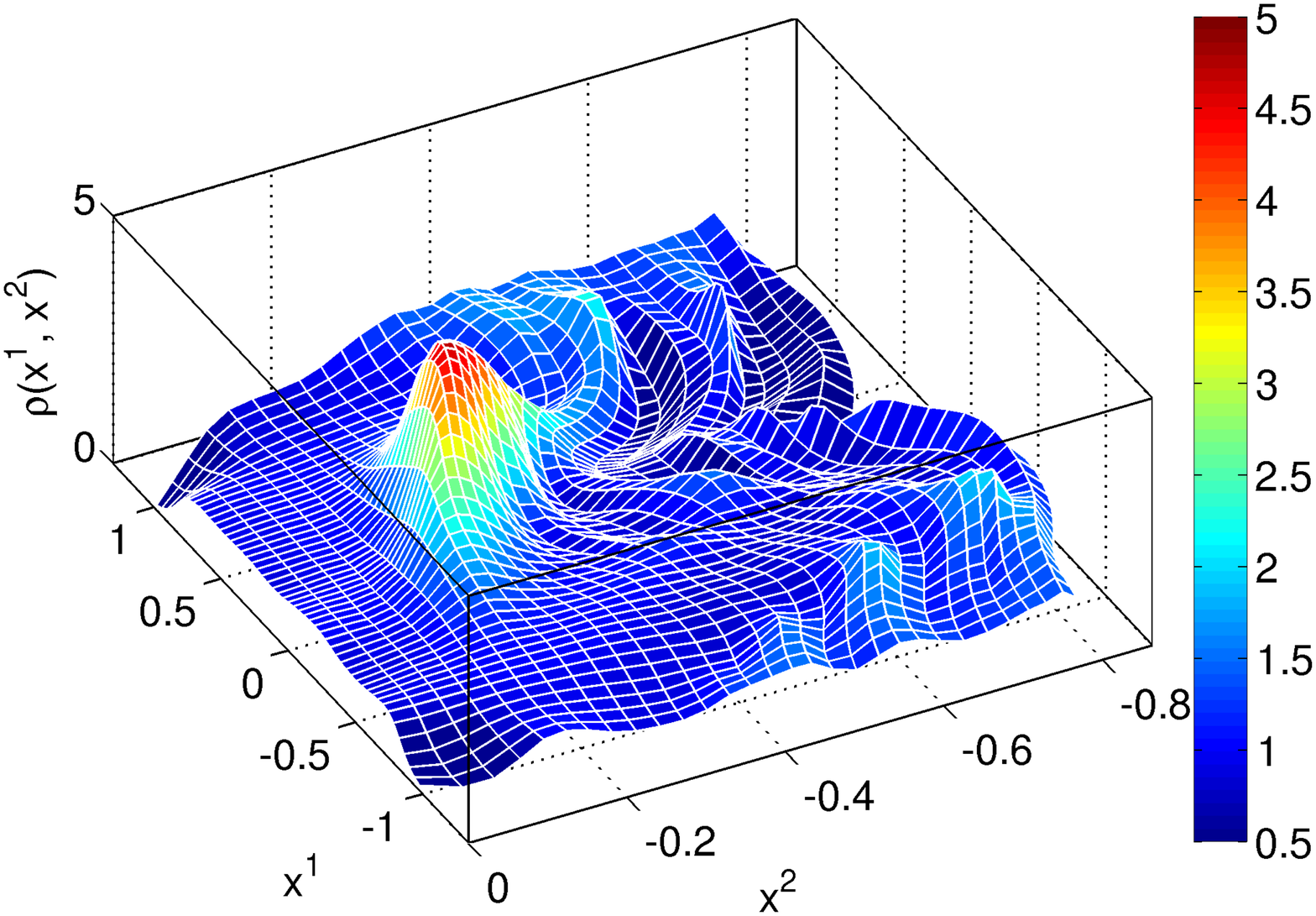}
\caption[]{Test 5. Density of medium $\rho(x)$ in the domain filled by waves initiated from $\sigma$:
left - exact values, right - recovered values.}
\label{fig:case5:dens}
\end{center}
\end{figure}

\noindent{\bf Test 5.}\,\,\, The purpose of the last test is to
check the ability of BCM to work with strong gradients in the
recovered quantities. We prepare the density of medium as a
slightly smoothed wedge with density $\rho = 5$ included in the
homogeneous background with the constant density $\rho = 1$, see
Figure \ref{fig:case5:dens}.

We take $T = 1$ and use the basis with 16 spatial and 32 temporal
functions. The condition number of matrix (\ref{Gramm G1}) for
$\xi = T$ is $2.07{\cdot}10^4$ and in calculations the parameter
of Tikhonov regularization for all linear systems is
$1.0{\cdot}10^{-4}$. The standard deviations of the Gaussian
kernels in (\ref{convol}) for $(\gamma, t)$ are $\sigma_\gamma =
0.125$ and $\sigma_t = 3.125{\cdot}10^{-2}$. Again, for
$\xi\approx T$ we had to increase $\sigma_\gamma$ to the value
$0.325$ in order to decrease large scale fake oscillations from
the low spatial harmonics.

\begin{figure}[ht]
\begin{center}
\includegraphics[width=0.54\textwidth]{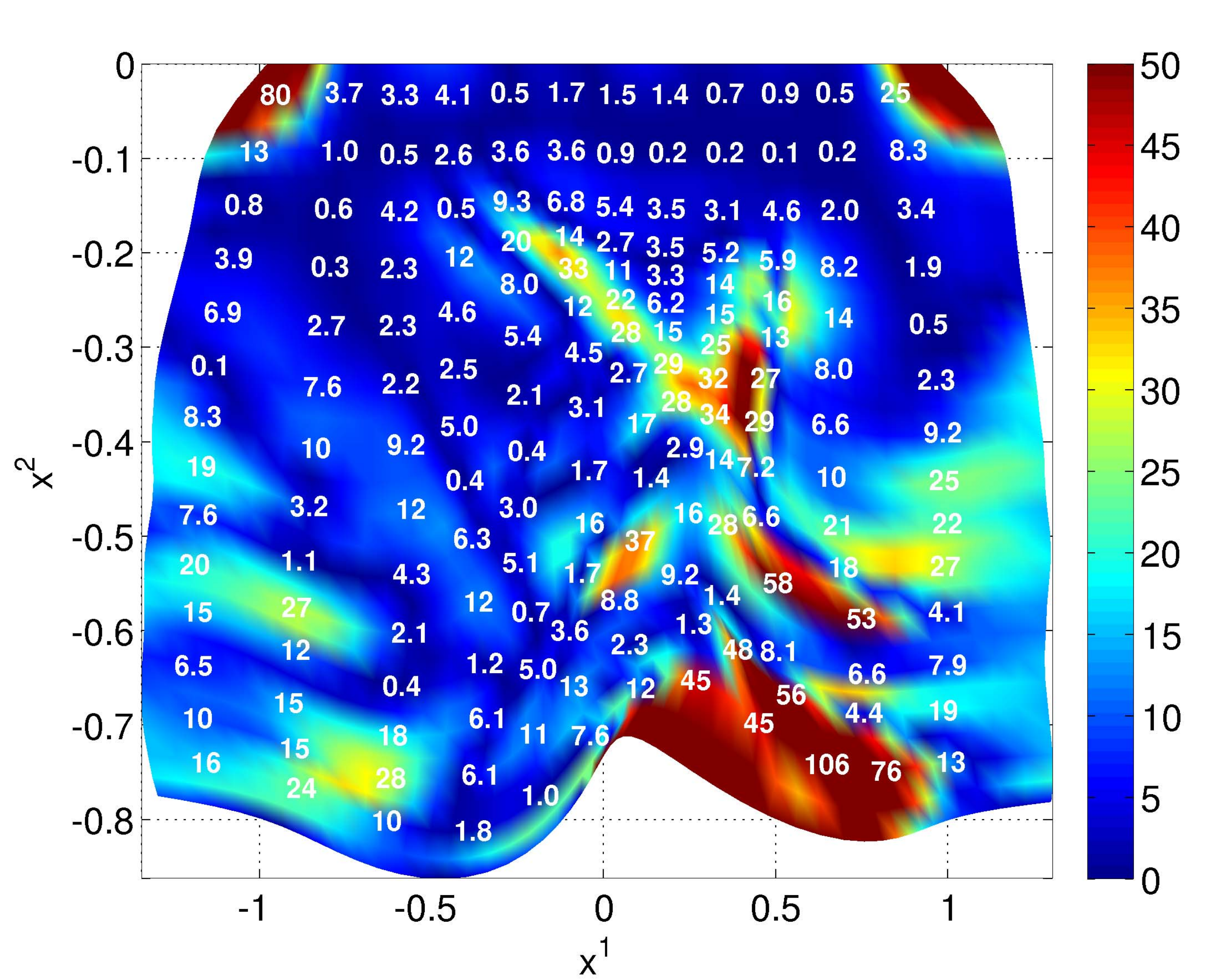}
\caption[]{Test 5. Map of relative errors of the recovered sound speed in percents.}
\label{fig:case5:relerr}
\end{center}
\end{figure}

The recovered density in the Cartesian coordinates is shown in
Figure \ref{fig:case5:dens}, and relative errors of reconstruction
for $c(x)$ in percents are shown in Figure \ref{fig:case5:relerr}.
The location of the wedge is recovered with good accuracy and
without systematic shifts, the maximum of recovered density is
$4.68$ that corresponds to the relative error about $7\%$. The
smearing of discontinuities is quite reasonable taking into
account the spatial scales of boundary controls and the smoothing
ranges, whereas the errors at the discontinuities are big, as it
has to be. We may suggest that if the resolution of boundary
controls is not enough to resolve spatial scales of the medium
inhomogeneities, then the method will recover a smoothed averaged
profile, which can be further improved by another high resolution
methods.

\subsubsection*{Comments}
$\bullet$\,\,\,Our results on the numerical speed determination
from the {\it time domain} data given at a {\it part} of the
boundary demonstrate that the BCM-algorithm is workable and
provides good reconstruction in the domain covered by the normal
acoustic rays.
\smallskip

\noindent$\bullet$\,\,\,The key step of the algorithm is inversion
of a big-size Gram matrix ($N\sim 10^2-10^3$), which consists of the inner products of
waves initiated by rich enough system of boundary controls. As is
typical in multidimensional (strongly ill-posed) inverse problems,
the condition number of this matrix rapidly grows as one extends
the number of controls and/or the observation time. Therefore, to
increase the depth of reconstruction one has to use controls
acting from a larger part of the boundary or decrease errors in
the input inverse data.
\smallskip

\noindent$\bullet$\,\,\,The number and shape of boundary controls
determine the spatial resolution of the reconstruction procedure.
The BCM is able to work with low number of spatial controls: in
such a case it provides an `averaged' profile. As we hope, such a
profile can be used as a starting approximation for high
resolution iterative reconstruction methods \footnote{Such an
option was suggested to the authors by F.Natterer.}.

In the future work, we plan to evaluate the influence of external
noises in the inverse data on the quality of the BCM
reconstruction, and apply the method to another domains and more
realistic sound speed profiles.

\newpage


\begin{thebibliography}{99}
\bibitem{ABI}
S.A.Avdonin, M.I.Belishev, S.A.Ivanov.
\newblock {Controllability in filled domains for the wave equation with singular controls.}
\newblock {\em Zapiski Nauchnykh Seminarov POMI},
210 (1994), 3--14 (in Russian); English translation: JMS,
83(1997), no 2.

\bibitem{BB}
V.M.Babich and V.S.Buldyrev,
\newblock {Short-Wavelength Diffraction Theory. Asymptotic Methods.}
\newblock {\em Sprin\-ger-Ver\-lag, Berlin}, 1991.

\bibitem{BeiKl},
L.Beilina and M.V.Klibanov.
\newblock {Approximate Global Convergence and Adaptivity for
Coefficient Inverse Problems.}
\newblock {\em Springer}, 2012.

\bibitem{BDAN87}
M.I.Belishev.
\newblock {On an approach to multidimensional inverse problems for
the wave equation.}
\newblock {\em Soviet Mathematics. Doklady}, 36 (1988), no 3, 481--484.

\bibitem{BIP97}
M.I.Belishev.
\newblock {Boundary control in reconstruction of manifolds and
metrics (the BC-method).}
\newblock {\em Inverse Problems}, 13 (1997), no 5, R1--R45.

\bibitem{BHow02}
M.I.Belishev.
\newblock {How to see waves under the Earth surface (the BC-method for
geophysicists).}
\newblock {\em Ill-Posed and Inverse Problems},
\newblock {S.I.Kabanikhin and V.G.Romanov (Eds). VSP, Utrecht, Boston}, 67--84, 2002.

\bibitem{BIP07}
M.I.Belishev.
\newblock {Recent progress in the BC-method.}
\newblock {\em Inverse Problems}, 23 (2007), no 5, R1--R67.

\bibitem{CUBO}
M.I.Belishev.
\newblock {Dynamical Inverse Problem for the Equation $u_{tt}-\Delta u -
\nabla \rho \cdot \nabla u=0$ (the BC-method).}
\newblock {\em CUBO A Mathematical Journal}, 10 (2008), No 2, 17--33.

\bibitem{BBlag99}
M.I.Belishev, A.S.Blagoveschenskii.
\newblock {Dynamical Inverse Problems of Wave Theory.}
\newblock {SPb State University, St-Petersburg}, 1999
(in Russian).

\bibitem{BGotJIIPP}
M.I.Belishev, V.Yu.Gotlib.
\newblock {Dynamical variant of the BC-method: theory and numerical testing.}
\newblock {\em Journal of Inverse and Ill-Posed Problems}, 7,
no 3: 221--240, 1999.

\bibitem{BGICOCV}
M.I.Belishev, V.Yu.Gotlib, S.A.Ivanov.
\newblock {The BC-method in multidimensional spectral inverse problem:
theory and numerical illustrations.}
\newblock {\em Control, Optimization and Calculus of Variations}, 2:
307--327, October, 1997.

\bibitem{BRF}
M.I.Belishev, V.A.Ryzhov, V.B.Filippov.
\newblock {Spectral variant of the BC-method: theory and numerical experiment.}
\newblock {\em Doklady Akad. Nauk SSSR}, 332 (1994), No 4,
414--417 (in Russian). English translation: ?????.

\bibitem{Ik}
M.Ikawa.
\newblock {Hyperbolic PDEs and Wave Phenomena.}
\newblock{\em Translations of Mathematical Monographs, v. 189}
\newblock{AMS; Providence. Rhode Island}, 1997.

\bibitem{Isbook}
V.Isakov.
\newblock{Inverse problems for partial differential equations.}
\newblock{\em Appl. Math. Studies, Springer}, v. 127, 1998.

\bibitem{KSh}
S.I.Kabanikhin, M.A.Shishlenin, A.D.Satybaev.
\newblock{Direct Methods of Solving Inverse Hyperbolic Problems.}
\newblock{\em Utrecht, The Netherlands, VSP}, 2004.

\bibitem{KSh_1}
S.I.Kabanikhin and M.A.Shishlenin.
\newblock{Numerical algorithm for two-dimensional inverse
acoustic problem based on Gel'fand-Levitan-Krein equation.}
\newblock {\em Journal of Inverse and Ill-Posed Problems}, 18,
no 9: 221--240, 2011.

\bibitem{kurganov2001}
A.Kurganov, S.Noelle, and G.Petrova.
\newblock{Semidiscrete central-upwind
schemes for hyperbolic conservation laws and Hamilton--Jacobi
equations.}
\newblock{\em SIAM Journal on Scientific Computing}, 23 (2001), no
3, 707--740.

\bibitem{Oks}
L.Oksanen.
\newblock {Solving an inverse obstacle problem for the wave equation by using
the boundary control method.}
\newblock {\em Inverse Problems}, 29 (2013), no 3, 035004; doi:10.1088/0266-5611/29/3/035004
http://iopscience.iop.org/0266-5611/29/3/035004/article?fromSearchPage=true.

\bibitem{P-12}
L.N.Pestov.
\newblock {Inverse problem of determining absorption
coefficient in the wave equation by BC method.}
\newblock {\em Journal of Inverse and Ill-Posed Problems}, 20,
no 1: 103--110, 2012.  ISSN (Online) 1569-3945, ISSN (Print)
0928-0219, DOI: 10.1515/jip-2011-0015, March 2012

\bibitem{P-13}
L.N.Pestov.
\newblock {On determining an absorption coefficient and a speed
of sound in the wave equation by the BC-method.}
\newblock {\em Journal of Inverse and Ill-Posed Problems}, 21,
no 2: 245--250, 2013. ISSN (Online) 1569-3945, ISSN (Print)
0928-0219, DOI: 10.1515/jip-2013-0012, November 2013.

\bibitem{PKBolg}
Pestov L., Kazarina O. Bolgova V.
\newblock {Numerical recovering a density by the boundary control method.}
\newblock {\em Inverse Problems and Imaging}, 2011. Vol. 4,
no. 4, 703--712.

\bibitem{Rom}
V.G.Romanov.
\newblock {A local version of a numerical
method for solving an inverse problem.}
\newblock {\em Siberian Math. J.}, 37 (1996), no 4, 797--810.

\end{thebibliography}
\end{document}